\documentclass[twoside,reqno,11pt]{amsart}
\usepackage{times,amsfonts,enumerate,amsmath,amsthm,amssymb,varioref}
\usepackage{hyperref}
\usepackage[alphabetic]{amsrefs}
\usepackage{array}
\usepackage{xy} \xyoption{all}
\usepackage[letterpaper,margin=4cm]{geometry}
{\end{list}}%

\newcount\hours
\newcount\minutes
\def\SetTime{\hours=\time
        \global\divide\hours by 60
        \minutes=\hours
        \multiply\minutes by 60
        \advance\minutes by-\time
        \global\multiply\minutes by-1 }
\SetTime
\def\now{\number\hours:\ifnum\minutes<10 0\fi\number\minutes}

\edef\today{\number\day\space
\ifcase\month\or
January\or February\or March\or April\or May\or June\or July\or August\or September\or October\or November\or December\fi, \number\year}

\def\Now{\today\space at \now}

\newcommand{\cancelqed}{\renewcommand{\qedsymbol}{}}

\def\addnotation #1: #2 #3{$#1$\> \parbox{5in}{\ \ #2 \dotfill\ \pageref{#3}}\\}

\providecommand{\gul}[1]{\protect{\langle#1\rangle}}
\providecommand{\norm}[1]{\protect{\lVert#1\rVert}}
\providecommand{\bignorm}[1]{\protect{\Big\lVert#1\Big\rVert}}
 \providecommand{\relnorm}[1]{\protect{\left\|#1\right\|}}

\providecommand{\bigabs}[1]{\protect{\big\lvert#1\big\lvert}}
\providecommand{\Bigabs}[1]{\protect{\Big\lvert#1\Big\lvert}}

\providecommand{\abs}[1]{\protect{\lvert#1\rvert}}

 \providecommand{\relabs}[1]{\protect{\left\lvert#1\right\lvert}}
 \providecommand{\relbrace}[1]{\protect{\left\{ #1 \right\}}}
 \providecommand{\relpar}[1]{\protect{\left( #1 \right)}}

\theoremstyle{plain} 
\newtheorem{teo}  {Theorem} 
\newtheorem*{teorema.intro}  {Theorem}   
\newtheorem{propo}[teo]{Proposition} 
\newtheorem*{prop.intro}{Proposition}   
\newtheorem*{lemma.intro}{Lemma}   
\newtheorem{lemma}[teo]{Lemma}       
\newtheorem{coro}[teo] {Corollary}  

\theoremstyle{definition}
\newtheorem{remark} [teo]  {Remark} 
\newtheorem{defi}[teo]   {Definition}
\newtheorem*{defi.sin}  {Definition}   

\theoremstyle{remark}

\newcommand{\Op}[1]{\, \mathrm{Op}\left(#1\right)}
\newcommand{\Opc}[1]{\, \mathrm{Op_c}\left(#1\right)}
\newcommand{\supp}[1]{\, \mathrm{supp}\,#1}

\providecommand{\dd}[2] { \frac{\partial}{\partial #1_{#2} }}

\newcommand{\sombrero}{\protect{\widehat{\hspace*{2.5mm}}}}
\newcommand{\spcheck}{^\vee}
\newcommand{\copete}[1]{\check{#1}}
\newcommand{\bm}[1]{\mbox{\boldmath{$#1$}}}

\newcommand{\ecal}{\mathcal{E}}

\newcommand{\kcal}{\mathcal{K}}

\newcommand{\mcal}{\mathcal{M}}

\newcommand{\ep}{\ep}
\newcommand{\mcalcop}{\check{\mathcal{M}}}
\providecommand{\mcalcopp}[1]  {\mcalcop^{\mathstrut^{#1}}}
\providecommand{\peq}[1]  {{{\scriptscriptstyle\mathstrut^{#1}}}}
\providecommand{\subpeq}[1]  {{{\scriptscriptstyle\mathstrut_{#1}}}}

\providecommand{\mcall}[1]     {\mcal^{\mathstrut^{#1}}}

\newcommand{\scal}{{\mathcal{S}}}

\newcommand{\rn} { \protect{{\mathbb R}^n}      }
\newcommand{\rbb}{\protect{\mathbb R} }
\newcommand{\rr}{\protect{\mathbb R}}

\newcommand{\z}{\protect{\mathbb Z}}

\newcommand{\frakg}{\protect{\mathfrak{g}}}

\newcommand{\frakX}{\protect{\mathfrak{X}}}
\newcommand{\frakA}{\protect{\mathfrak{A}}}
\newcommand{\frakB}{\protect{\mathfrak{B}}}

\providecommand{\sro}[1]  {S^m_{\rho,{#1}}}
\providecommand{\sroo}[2]  {S^{#1}_{\rho,{#2}}}
\newcommand{\cisc}{C^\infty_c }

\newcommand{\ci}{C^\infty}

\newcommand{\id}{\mathop{\mathrm{Id}}}

\numberwithin{equation}{section}

\newlength{\anchocajateo}
\numberwithin{equation}{section}
\begin{document}

 \date{\Now.}

\author{Susana Cor\'e}
\email{score@math.smith.edu}
\address{Department of Mathematics\\
Smith College\\ 
  Northampton, MA 01063}

\author{Daryl Geller }
\email{daryl@math.sunysb.edu}
\address{Department of Mathematics, 
          Stony Brook University,
          Stony Brook, NY 11794-3651}


\title[Pseudodifferential Calculus on Homogeneous Groups]{H\"ormander Type Pseudodifferential Calculus on Homogeneous Groups}

\begin{abstract}
 We produce, on general homogeneous groups, an analogue 
 of the usual H\"ormander pseudodifferential calculus 
  on Euclidean space, at least as far as products and
  adjoints are concerned.  In contrast to earlier works, 
  we do not limit ourselves to analogues of classical symbols, 
   nor to the Heisenberg group.  
   The key technique is to understand \emph{``multipliers''} 
  of any given order $j$, and the operators of convolution 
   with their inverse Fourier transforms, which we here
  call convolution operators of order $j$.  
  (Here a \emph{``multiplier''} is an analogue of a 
  H\"ormander-type symbol $a(x,\xi)$, which is independent 
  of $x$.)
  Specifically, we characterize the space of inverse Fourier 
  transforms of multipliers of any order $j$, and use 
  this characterization to show that the composition of 
   convolution operators of order $j_1$ and $j_2$ is 
   a convolution operator of order $j_1+j_2$.
\end{abstract}

\maketitle
\setcounter{tocdepth}{1}
\tableofcontents


\section{\sc Introduction}
The goal of this article is to produce, on general homogeneous groups 
  (in the sense of Folland-Stein\cite{folland-stein-libro}),
  an analogue of the {\em usual}  pseudodifferential calculus 
  on Euclidean space, at least as far as products
  and adjoints are concerned.  Taylor\cite{taylor-part-i} and
  Christ-Geller-G\l owacki-Polin \cite{geller-pollin}
  have developed analogues of \emph{classical}  pseudodifferential operators
  for homogeneous groups, but, as we shall explain, our analogue of 
  the usual calculus allows one to study a useful wider class of operators. 
 
The first author \cite{susana-heisenberg} has already produced 
an analogue of the {\em usual} pseudodifferential calculus 
for the Heisenberg group.  The arguments in the present article
are simpler even in this special case.  
(Note, however, that in \cite{susana-heisenberg}, the first author completely 
characterized the kernels of the
operators if they had negative order, 
even on general homogeneous groups $G\neq\rr$.  
This information on the kernels of the operators 
will be valuable in sequel articles.)\\

Before describing the methods we use, we make more precise 
the terms \emph{usual} and \emph{classical}
pseudodifferential calculus used loosely in the previous paragraphs.
We recall that 
 the {\em (H\"or\-man\-der) symbol class of order} $m$,
 denoted by $S^m_{1,0}$,
 consists of those functions $a(x,\xi)$ in $\ci(\rn\times\rn)$
 such that for any pair of multiindices $\alpha, \beta$, and any
 compact set $B\subset \rn$, there exists a constant
 $C_{\alpha,\beta,B}$ such that
     \begin{equation}\label{regsym}
   \Bigabs{D^\beta_x D_\xi^\alpha a(x,\xi)}
   \leq
   C_{\alpha,\beta,B} \big( 1+\norm{\xi} \big)^{m - \norm{\alpha}}
\qquad \forall x\in B,\ \xi\in\rn
    \end{equation}
Given that H\"or\-man\-der's $S_{1,0}$ symbols have become standard 
we shall refer to them as the  {\em usual symbols,}
 and the associated calculus as the {\em usual} or
{\em H\"or\-man\-der calculus} of type $(1,0)$.
 The {\em classical symbols of order} $m$
 are those elements $a(x,\xi)$ in $S^m_{1,0}$  
  for which there are smooth functions 
 $a_{m-j} (x,\xi)$, homogeneous of degree $m-j$
in $\xi$,  for $\norm{\xi} \geq 1,$ 
 such that
 \begin{equation*}\label{}
    a(x,\xi) \sim  \sum_{j\geq 0} a_\subpeq{m-j} (x,\xi)
    \end{equation*}
 where the asymptotic condition means that for any $N$
\begin{equation*}\label{}
   a(x,\xi) -\sum_{j=0}^N a_{m-j} (x,\xi) \in S^\peq{m-N-1}_{1,0}
 \end{equation*}
  Our approach to pseudodifferential operators on homogeneous groups favors
 the use of  convolution operators
and avoids the Fourier transform as much as possible.  This idea 
goes back  to the original work of Mikhlin \cites{mikhlin,mikhlin-trad} 
and of Cal\-de\-r\'on and Zygmund 
 \cite{calderon-zygmund}, \cite{calderon-zygmund-dos}.
    These pioneers studied operators of the form 
\begin{equation*}
     (\kcal f) (x) = \int K(x, x-y) \, f(y) \, dy
\end{equation*}
so that
\begin{equation}\label{calde}
  (\kcal f) (x) = ( K_x \ast f ) (x)
    \end{equation}
with $K_x(z) = K (x,z)$ 
an integral kernel for each $x$, 
singular only at the origin, and depending smoothly on $x$.
In the 1960's Kohn and Nirenberg \cite{kohn-nirenberg} 
introduced symbols, a different point of view,  
and the term {\em pseudodifferential operator}. 
They used the Fourier transform to rewrite 
the  operator $\kcal$ in the form
\begin{equation*}\label{kohn}
         (\kcal f) (x)
   = \int  e^{-2\pi i x \cdot \xi}\ a(x,\xi) \ \widehat{f}(\xi) \, d\xi
   = \left( a_x \widehat{f} \right)\check{\ \ } (x)
  \end{equation*}
where $a$, the  symbol of the operator,
is the formal Fourier transform of $K$ in the second variable, that is
\begin{equation}\label{axkxhat}
 {a}_x ( \xi ) = \widehat{K}_x(\xi ) \quad \text{and}\quad  a_x(\xi) = a(x,\xi)
    \end{equation}
This last approach, which takes advantage of the properties of the Fourier
transform, became prevalent.  A major reason is that the Fourier transform converts
convolution to a product, which is easier to handle.  In particular,
one can seek to use division to invert $\kcal$.

In problems where the Euclidean convolution structure is not relevant,
these advantages are largely lost, and in many instances it is
 desirable
to imitate the original definition of Mikhlin and Calder\'on-Zygmund.

For instance,
if one is working on a Lie group, one can seek to define a class of
pseudodifferential operators by (\ref{calde}), where now $*$
is group convolution.  This idea originates in Folland-Stein 
\cite{folland-stein-estimates}
for the Heisenberg group, and was extended in Rothschild-Stein
 \cite{rothschield-stein} to other settings.  Calculi of such operators 
on homogeneous groups were
developed by Taylor \cite{taylor-part-i} and studied in greater 
detail in \cite{geller-pollin}.  (A related calculus, restricted to the Heisenberg
group, and relying heavily on the Fourier transform, was developed in
\cite{beals}.)  The second author further developed these ideas,
to obtain a calculus in the real analytic setting on the Heisenberg group, in 
\cite{geller-libro}.   All of these authors restricted themselves to analogues
of the classical calculus.

The operators in the calculus of this article have the form 
(\ref{calde}) (with $*$ being group convolution), 
with $K_x$ being again given by (\ref{axkxhat}),
but where now $a(x,\xi)$ satisfies, for any compact set $B \subseteq \rn$, 
the estimates
     \begin{equation}\label{gpsym}
   \Bigabs{D^\beta_x D_\xi^\alpha a(x,\xi)}
   \leq
   C_{\alpha,\beta,B} \big( 1+|\xi|)^{m - |\alpha|}
\qquad \forall x\in B,\ \xi\in\rn
    \end{equation}
in analogy to (\ref{regsym}).  (Let us then say that $\kcal$ has order $m$.)
It is crucial to note that in this definition, and in what follows, 
$\abs{\xi}$ denotes the homogeneous norm of $\xi$, and $\abs{\alpha}$
is the ``weighted'' length of the multiindex $\alpha$.  

Thus our 
$a(x,\xi)$ are not required to have an asymptotic expansion in 
homogeneous functions of $\xi$, as $\xi \rightarrow \infty$, as is required
in the classical calculi.  

We will show that there is a calculus for these
operators, in that the adjoint of an operator of order $m$ is again 
an operator of order $m$, and that the composition of operators of 
order $m_1$ and $m_2$ is an operator of order $m_1+m_2$.  (When composing,
one assumes that the ``$a$'' associated to the operator which is applied first, has
compact support in $x$.)  Further, there are asymptotic expansions, similar to the
Kohn-Nirenberg formulas, for adjoints and compositions.
 
\subsection{\sc outline and list of results}
This paper is organized as follows: 
\begin{enumerate}[\phantom{mm}  $\bullet$]
\item
Section 1 contains this introduction, and Section 2 is dedicated to establishing notation and basic terminology.
\item
   In Section 3 
   we study the spaces $\mcalcop^j(G)$, consisting of inverse 
    Fourier transforms of {\em ``multipliers''}.  A ``multiplier'',
   by definition, is a function $a$ as in (\ref{gpsym}) which is independent 
   of $x$.  Thus:
\begin{defi.sin}
 Suppose $m\in \rr$. We shall say that $u\in\ci(G)$ is
  a \emph{ multiplier of order $m$}  if for every
  multiindex $\alpha\in (\z^+)^n$   there exists
  $C_\alpha >0$  such that
      $$
       \bigabs{\partial^\alpha u (\xi)}
        \leq
       C_\alpha (1+\abs{\xi})^{m-\abs{\alpha}}
      \qquad \qquad  \text{ for all } \xi
      $$
$\mcal^{m}(G)$  will denote the space of multipliers of order $m$.
 \end{defi.sin}

If $a$ satisfies (\ref{gpsym}), and $a_x(\xi) = a(x,\xi)$, then for any $x$, 
   we surely have $a_x \in {\mathcal M}^m$.  The operator $\kcal$ of
   (\ref{calde}), (\ref{axkxhat}), is then of the form 
   $\kcal f(x) = [\check{a}_x * f](x)$.

In this sense, our calculus of pseudodifferential operators 
   will have, as its building blocks, operators of
   the form $Af = \check{u}*f$, where $u \in \mcal^{j}(G)$ for some $j$
   and $*$ is convolution on $G$, and it is these convolution operators that
   will require our major focus.  To avoid confusion, we caution the reader 
   that these convolution operators will only be Fourier multiplier operators 
   (in the ordinary sense) if $*$ is Euclidean convolution.  
   We nevertheless call elements of $\mcal^{j}(G)$ 
   \emph{``multipliers'',} for lack of a better word.  

The main result of Section 3 is the following characterization of 
   the space of the inverse Fourier transform of multipliers, 
  $\mcalcopp{j}(G)$.  Here $Q$ denotes the homogeneous dimension 
  of $G$, $\scal(G)$ denotes the space of Schwartz functions 
  on $G$ (which coincides with the usual space of Schwartz 
   functions on $\rn$), and $\scal_o(G)$ denotes the space 
  of elements of $\scal(G)$, all of whose moments vanish.
\begin{prop.intro}
   Say $K \in \scal'(G)$.  Then $K \in \mcalcopp{j}(G)$ 
   if and only if we may write
   \begin{equation*} \label{kmgcrt.pre}
     K = \sum_{k=0}^\infty f_k,
   \end{equation*}
   with convergence in $\scal'$, where
   \begin{equation} \label{fkphi.pre}
   f_k (x) = 2^{(j+Q)k} \varphi_k (\delta_{2^k}(x))
   \end{equation}
   where    $\left\{\varphi_k \right\}_k \subseteq \scal(G)$  
   is a bounded sequence, and 
$\varphi_k \in \scal_o(G)$ for $k \geq 1$.    Further, given 
   such a sequence $\left\{\varphi_k \right\}_k$, 
   if we define $f_k$ by (\ref{fkphi.pre}), then 
   the series $\sum_{k=0}^\infty f_k$ necessarily
   converges in $\scal'(G)$ to an element of $K$ of $\check{\mcal}^j (G)$.
\end{prop.intro}

Related decompositions (with far weaker conditions on the $\varphi_k$) 
   have occurred before in the literature on multipliers, when minimal 
  smoothness of the multiplier was assumed.  (In the situation of general 
  dilations, but where ordinary Euclidean convolution is used, 
  see \cite{ricci-multipliers}, pages 37-43.) It is our intention to show that 
  this kind of decomposition is extremely helpful even when one assumes 
 the full force of the $\mcal^j$ conditions.
\item
  In Section 4 the main result is the following key theorem:
\begin{teorema.intro}
 For any $j_1, j_2$,   
          $$\mcalcop^{j_1}(G)\ast \mcalcop^{j_2}(G) 
            \subseteq
            \mcalcop^{j_1+j_2}(G).
             $$
\end{teorema.intro}
The proof uses the characterization of $\mcalcop^j(G)$ given in Section 3, an adaptation of Lemma 3.3
of Frazier-Jawerth \cite{frazier-jawerth} to the homogeneous group setting, 
and some additional ideas.  

Here is our lemma, adapted from Frazier-Jawerth \cite{frazier-jawerth}, 
who established the case in which $G$ is Euclidean space.
For $\sigma \in \z$, $I > 0$, define
\[ \Phi^I_{\sigma}(x) = (1+2^{\sigma}|x|)^{-I}. \]
\begin{lemma.intro}
Say $J > 0$, and let $I = J+Q$.  Then there exists $C > 0$ such that whenever $\sigma \geq \nu$, 
\[ 
   \Phi^I_{\sigma}*\Phi^J_{\nu} 
   \leq 
    C\, 2^{-\sigma Q} \Phi^J_{\nu}. 
\]
\end{lemma.intro}
Again, Frazier and Jawerth proved their lemma in order to study spaces 
   of restricted smoothness (specifically, Besov spaces).  
   But we intend to show that this lemma is also useful 
   when one assumes the full force of the $\mcal^j$ conditions.
\item
 In Section 5, we obtain our new calculus.  At this point we simply
refer to the general theory of Taylor\cite{taylor-part-i} to
say that the desired calculus exists, as long as one has the
key result of Section 4 (that $\mcalcop^{j_1}(G)\ast \mcalcop^{j_2}(G) \subseteq
\mcalcop^{j_1+j_2}(G)$), together with a few easily checked facts
(which we do verify).
\end{enumerate}

In future articles, we will examine the calculus in more detail.
  Among other properties, we will seek to obtain, as in the classical 
  situation discussed in detail in \cite{geller-pollin}, explicit formulas 
 for products and adjoints, criteria for existence of parametrices, 
  and mapping properties of our pseudodifferential operators.

As we said, all of those properties are known for the analogues 
   of {\em classical} pseudodifferential operators on $G$.  
   Let us then motivate our present work by mentioning two expected 
  applications of our new calculus, that the classical calculus 
  cannot be expected to have.  
\begin{enumerate}[\phantom{mm}  $\bullet$]
\item
Classical multipliers are, as $\xi \rightarrow \infty$, 
   asymptotic sums of series of terms which are homogeneous 
   with respect to the (non-isotropic) dilations on $G$.  
   Let us call these dilation $\delta_r$ (for $r > 0$).  
   One does not always want to restrict oneself to functions 
   which are homogeneous with respect to $\delta_r$ for {\em all } $r > 0$; 
  in the theory of wavelets, for instance, one really only cares 
   about {\em dyadic} dilations.  
   It is easy to see that if a smooth function $u(\xi)$ 
    is homogeneous of degree $j$ with respect 
    to $\delta_2$ (for $\xi$ outside a compact set), 
    then it is in $\mcal^{j}(G)$.  

In \cite{geller-mayeli}, wavelet frames 
(for $L^p$, $1 < p < \infty$, and $H^1$)
   were constructed on stratified Lie groups $G$ with
   lattice subgroups, out of Schwartz functions $\psi$ of the form
   $f({\mathcal L})\delta$, for  a nonzero $f \in {\mathcal S}(\rr^+)$ 
   with $f(0) = 0$; here ${\mathcal L}$ is the sublaplacian on $G$.  
   Part of the argument involved inverting a spectral multiplier 
   (for the sublaplacian) which was homogeneous of degree zero 
   with respect to dyadic dilations, and using the fact that 
   it therefore satisfied standard multiplier conditions.   
   We expect that, once we understand inversion for convolution 
  operators of the kind considered in this article, we will be able 
  to generalize these arguments to much more general Schwartz 
   functions $\psi$ on $G$ (not just those of the special form
  $f({\mathcal L})\delta$), and therefore be able to construct 
   wavelet frames from much more general $\psi$.

\item
Let $(\bm{M},g)$ be a smooth, compact, oriented Riemannian manifold.  
   One of the most important facts about analysis on $\bm{M}$ is 
   the following theorem of Strichartz \cite{strichartz}:
   \begin{quote}
    \begin{center}
   \emph{
    If $ p(\xi) \in S^{m}_{1,\#}(\rr) $  
    then $p(\sqrt{\Delta}) \in OPS^m_{1,0}(\bm{M})$.
   }
    \end{center}
   \end{quote}
   ($S^m_{1,\#}$ denotes the space of symbols $a(x,\xi)$ 
   in $S^{m}_{1,0}$ which are independent of $x$,
   or in other words, multipliers of order $m$ in the usual sense.)  
   Here $\Delta$ denotes the Laplace-Beltrami operator on $\bm {M}$, 
   but one can replace $\sqrt{\Delta}$ by a general first order 
   positive elliptic pseudodifferential operator on $\bm{M}$.
  As in the classical case \cite{taylor-part-i}, \cite{geller-pollin}, 
   taking $G$ to be the Heisenberg group, we expect to be able 
   to transplant our calculus to contact manifolds, and in particular, 
   to smooth, compact strictly pseudoconvex CR manifolds $\bm{M_0}$.  
   We expect that one can show:
 \begin{quote}
 \emph{
   If $ p(\xi) \in S^{m}_{1,\#}(\rr)$,  then $ p(\sqrt{L}) $
    is an operator of order $ m$ 
     in the new calculus on  $\bm{M_0}$. 
   }
 \end{quote}
   Here $L$ denotes the sublaplacian on $\bm{M_0}$, and we expect 
   to be able to replace $\sqrt{L}$ here by a general first-order operator 
  in the new calculus on $\bm{M_0}$,
   which is transversally elliptic and positive. 
\ \\
\end{enumerate}

\section{\sc Preliminaries}
We present basic results for homogeneous groups, and introduce
the notation to be used later. 
For more details see~\cite{folland-stein-libro}.
\begin{defi}
  Let $V$ be a real vector space.  
  A family $\{\delta_r\}_{r>0}$ of linear maps of $V$
  to itself is called a set of {\em dilations on }$V$, 
  if there are real numbers $\lambda_j>0$ and subspaces
  $W_{\lambda_j}$ of V such that $V$ is the direct sum
  of the $W_{\lambda_j}$ and 
       $$
        {\delta_r}\big\vert_{W_{\lambda_j}}
         = r^{\lambda_j} \id \qquad
       \forall j 
        $$
\end{defi}
\begin{defi} 
 A {\em homogeneous group} is a connected and simply 
 connected nilpotent group $G$, with underlying manifold $\rn$,
 for some $n$, and 
 whose Lie algebra $\frakg$ 
is endowed with a family of dilations  $\{\delta_r\}_{r>0}$,
which are automorphisms of $\frakg$.
 \end{defi}
The dilations are of the form $\delta_r = \exp (A \log r)$,
where $A$ is a diagonalizable linear operator on $\frakg$
with positive eigenvalues.
The group automorphisms 
$ \exp \circ \delta_r \circ \exp^{-1}: G \longrightarrow G$
will be called {\em dilations} of the group and  
will also be denoted by $\delta_r$.
The group $G$ may be identified topologically with $\frakg$
via the exponential map $\exp: \frakg \longrightarrow G$
and with such an identification
 \begin{equation*}
   \begin{array}{rclc}
 \delta_r\, :& \, G              & \longrightarrow              & G\\
             &  (x_1,\ldots ,  x_n) & \longmapsto                  & 
          (r^{a_1} x_1, \ldots , r^{a_n} x_n)
   \end{array}
 \end{equation*}
Henceforth the eigenvalues of the matrix $A$, listed as many times as their
multiplicity, will always be denoted by $\big\{a_i\big\}_{i=1}^n$.
Moreover, 
we shall assume without loss of generality that 
 all the $a_i$ are nondecresingly ordered 
and that the first is equal to $1$, that is
   \begin{equation*}\label{} 
   1 = a_1 \leq \ldots  \leq a_n\, . 
    \end{equation*}
\begin{defi} 
The {\em homogeneous dimension} of the group $G$,
 denoted by $Q$, is the number   
\begin{equation*}\label{} 
   Q = \sum_{i=1}^n a_i \ .
 \end{equation*} 
 \end{defi}
Examples of such groups are $\rn$, with the usual additive
structure, the Heisenberg groups, 
and the upper triangular groups, consisting of 
triangular matrices with 1's on the diagonal 
and dilations
   \begin{equation*}\label{} 
    \delta_r \Big( \big[a_{ij}\big]\Big) 
    = 
    \big[r^{j-i} a_{ij}\big]
    \end{equation*}
 \begin{defi} 
   Let $G$ be a homogeneous group with dilations  $\big\{ \delta_r \big\}$.
   A {\em homogeneous norm} on $G$, relative to the given dilations, is
   a continuous function $\abs{\ \cdot \ } : G\longrightarrow [0,\infty)$,
   smooth away from the origin satisfying
   \begin{enumerate}[\ \ \ \  a)] 
  \item  
      $\abs{x} = 0$ if and only if $x$ is the identity element,
  \item $\abs{x^{-1}} = \abs{x}$ for every $x\in G$
  \item $\bigabs{\delta_r(x)} = r\abs{x}$ for every $x\in G$, and $r>0$,
   i.e. the norm is homogeneous of degree~$1$.
 \end{enumerate}
 \end{defi}
 If  $\abs{\ \cdot \ }$ is a homogeneous norm on $G$ then there exists
 a constant $C\geq 1$ such that 
  \begin{equation*}\label{} 
       \abs{x y } \leq C \big( \abs{x} + \abs{y} \big)
      \qquad \text{ for every } x, y \in G 
    \end{equation*}
 Homogeneous norms always exist.  Moreover any two homogeneous 
 norms  $\abs{\ \cdot \ }$ and  $\abs{\ \cdot \ }'$
 on $G$ are always equivalent, i.e. there exist
  constants $C_1, C_2> 0$ such that 
  \begin{equation*}\label{} 
       C_1 \abs{x} \leq \abs{x}' \leq C_2 \abs{x} \qquad \forall x\in G
    \end{equation*}
 An example of a homogeneous norm is 
  \begin{equation*}\label{} 
 \abs{x} = \left(  \sum_{i=1}^n x_i^\peq{2\frac{A}{a_i}}\right)^{\frac{1}{2A}} 
    \end{equation*}
    where  $A=\prod_{i=1}^n a_i$.

  We observe that if the dilations are isotropic, i.e. all the weights
  $a_i$ are equal, and satisfy the condition of normalization $a_1=1$,
  then this homogeneous norm is simply 
   the Euclidean norm $\norm{\, \cdot \, }. $

 If $x \in G$, $r>0$, 
 for a fixed homogeneous norm
 we define the \emph{ball centered at $x$ of radius $r$} as 
  \begin{equation*}\label{} 
   B(x;r) = \Big\{  
           y \in G \, : \, \abs{x^{-1} y} < r 
           \Big\} 
    \end{equation*}
Note also that there
exists $c > 0$ such that whenever $x,y \in G$, we have
\begin{equation}
    \label{yxest}
    |y| \leq |x|/2 \Rightarrow |y^{-1}x| \geq c|x|.
\end{equation}
Indeed, by use of dilations, we may assume that $|x| = 1$.  
We then let $c$ be the minimum value
of $|y^{-1}x|$ for $(x,y)$ 
in the compact set 
$\{x: |x| = 1\} \times \{y: |y| \leq 1/2\} \subseteq G \times G$.

For $j=1,\ldots,n$, we let   $X_j$ (resp.  $Y_j$)
denote the left (resp. right) invariant vector field  on $G$ 
which equals $\dd{x}{j}$ at the origin. In this context
we have
\begin{propo}\label{proposicion-pagina-41}
   \begin{enumerate}[a)] 
     \item 
     We may write each left invariant vector field $X_j$ as
         \begin{equation*}
         X_j 
         = 
         \dd{x}{j} + \sum^n_{k=j+1} p_{j,k} (x) \dd{x}{k}
         \qquad j=1,\ldots,n
         \end{equation*}
         where $p_{j,k} (x) = p_{j,k} (x_1,\ldots, x_{k-1}) $
         are homogeneous polynomials, with respect to the group dilations,
         of degree $a_k-a_j$.
     \item Any $\dd{x}{j}$ can be written as
         \begin{equation*}
         \dd{x}{j} = X_j + \sum^n_{k=j+1} q_{j,k} (x) X_k
         \qquad j=1,\ldots,n
         \end{equation*}
      where $q_{j,k}(x) = q_{j,k}(x_1,\ldots,x_{k-1}) $
      are homogeneous polynomials, with respect to the group dilations,
       of degree $a_k-a_j$.
      (Note:
        by homogeneity considerations $q_{j,k} (x)$ can only  
        involve $x_1, \ldots, x_{k-1},$ so
      by part a), multiplication by  $q_{j,k}(x)$
      commutes with $X_k$).
   \end{enumerate}
   Entirely analogous formulas express
   each  of the right invariant vector fields $Y_j$ 
   in terms of $\left\{\dd{x}{k} \right\}_{k=j}^n$,
   as well as
   $\dd{x}{j}$ in terms of $\left\{Y_{k}\right\}_{k=j}^n$
\end{propo}
\begin{proof}
 See \cite{folland-stein-libro}.
\end{proof}

The Haar measure on $G$ is simply the Lebesgue measure on $\rn$.
\begin{propo}
  Say  $p < -Q$. Then for some $0 < C_p < \infty $ 
   \begin{equation*}\label{} 
 \int_{\abs{x}>r} \bigabs{x}^p \ dx  =  C_p \, r^{p+Q} 
 \qquad 
    \end{equation*}
 for every $r > 0$.
 \end{propo}
 The \emph{convolution} of two functions $f,g$ on $G$ is defined
 by 
 \begin{equation*}\label{} 
   \big( f\ast g\big) (x) 
   = 
   \int_G f (xy^{-1}) g(y) \, dy 
   =  
   \int_G f(y) g (y^{-1}x)  \, dy 
    \end{equation*}
  provided that the integrals converge.
\noindent

$\scal(G)$ will denote the usual Schwartz space on $G$, 
thought of as $\rn$,
and $\big\{ \norm{\ \cdot\ }_{S(G), N}\big\} $ 
an increasing family of norms topologizing $\scal$.
We also let
$$ \scal_o(G) = \left\{f \in \scal(G): \int x^{\alpha} f(x) dx = 0 
  \mbox{ for every multiindex } \alpha\right\}. $$
Thus\footnote{
      It will be convenient to think of the Fourier
      transform on $\rn$ as a map from functions  on $G$
      to functions on $G$.
      To clarify this, on the Fourier transform side
      we will often use the dilations $\delta_r$  
      and the homogeneous norm $\abs{\ \cdot \ } $,
       but never the group structure.
        } 
$$ \widehat{\scal}_o (G) 
    = \left\{\psi \in \scal(G): \partial^{\alpha}\psi(0) = 0 
   \mbox{ for every multiindex } \alpha\right\}. 
$$
Since the group law is polynomial, $\scal_o$ is invariant
 under (left or right) translations. 
Thus
 $\scal \ast \scal_o \subset \scal_o$,
since if we write 
  $(\tau_y g)(x) = g (y^{-1} x) $,
  we have $f \ast g = \int f(y) \ \tau_y(g) \ dy$.
 Simlarly, $\scal_o  \ast \scal \subset \scal_o$

For the multiindex $\beta \in (\z^+)^n$
we define 
\begin{equation*}\label{} 
\abs{\beta}= \sum_{i=1}^n a_i \beta_i
\qquad\text{and}
\qquad
\norm{\beta}= \sum_{i=1}^n  \beta_i
    \end{equation*}
Note that, if $\alpha$ is a multiindex, there exists $C > 0$ such that
$|x^{\alpha}| \leq C|x|^{|\alpha|}$ for all $x \in G$. 
\begin{lemma}\label{lema-16}
 For every multiindex  $\alpha$
 \begin{equation*}\label{Eins} 
   a_1 \norm{\alpha} \leq \abs{\alpha} \leq a_n \norm{\alpha}
    \end{equation*}
  and for all $\xi\in \rn$ there exist positive constants $c$ and $C$
such that
  \begin{equation*}\label{Zwei} 
    c \big(1+\abs{\xi}\big) \leq \gul{\xi} \leq C \big(1+\abs{\xi}\big)^{a_n}
    \end{equation*}
 where $ \gul{\xi} =  \big(1+\abs{\xi}^2\big)^{1/2}$.
\end{lemma}


\begin{defi} \label{def-multiplicador}
 Suppose $j\in \rr$. We say that $J\in\ci(G)$ is
  a {\em multiplier of order $j$}  if for every 
  multiindex $\alpha\in (\z^+)^n$   there exists
  $C_\alpha >0$  such that
      $$ 
       \bigabs{\partial^\alpha J (\xi)}
        \leq
       C_\alpha (1+\abs{\xi})^{j-\abs{\alpha}}
      $$ 
  for all $\xi$. We denote the space of all multipliers of order $j$
  by $\mcall{j}(G)$.
 \end{defi}

  For $J\in \mcall{j}(G)$, and $N\in \z^+$, we define
  \begin{equation*}\label{} 
   \norm{J}_{\mcall{j}(G),N} = \sum_{\abs{\alpha}\leq N} 
                 \bignorm{ \big(1+\abs{\xi}\big)^{\abs{\alpha}-j}\
                               \partial^\alpha J
                         }
    \end{equation*}
  where $\norm{\ \cdot \ } $ denotes the supremum norm.
  The family $\left\{ \norm{\ \cdot \ }_{\mcall{j}(G),N} \right\}_{N\in\z^+}$
  is a nondecreasing
 sequence of 
 norms on $\mcall{j}(G)$,
 which defines a Fr\'echet topology on $\mcal^j(G)$.

 Suppose $J\in \mcal^j(G)$. 
 We define 
  $$ m_J: \scal(G)\longrightarrow \scal'(G) $$
  $$ m_J(f) = \check{J} \ast f$$

We close this section with some useful remarks
about $\scal(G)$ and $\scal_o(G)$.
\begin{remark}
\label{rmk1}
Note that for every $\psi \in \hat{\scal}_o$, 
  we may find $\psi_1,\ldots,\psi_n \in \hat{\scal}_o$ 
with $\psi(\xi) = \sum_{j=1}^n \xi_j \psi_j(\xi)$.  
  Indeed, we need only show this in two cases:
 (1) if $\psi$ vanishes in a neighborhood of $0$; 
 and
 (2) if $\psi \in C_c^{\infty}(G)$.  
 In case (1), it suffices to set 
    $\psi_j(\xi) = \xi_j \psi(\xi)/\|\xi\|^2$, 
 where $\norm{\xi } $ is the Euclidean norm of $\xi$.
In case (2), choose $\zeta \in C_c^{\infty}(G)$
with $\zeta \equiv 1$ in a neighborhood of $\supp \psi$.  
Fix $\xi$ temporarily, and let
$g(t) = \psi(t\xi)$.  Then
\[ 
     \psi(\xi) = g(1)-g(0) = \int_0^1 g'(t)dt 
      = \sum_{j=1}^n \xi_j \int_0^1 \partial_j \psi(t \xi) dt. \]
Noting that $\psi = \zeta \psi$, we see that in case (2) we may take 
$\psi_j(\xi) = \zeta(\xi)\int_0^1 \partial_j \psi(t \xi) dt$.

In fact, this construction shows the following uniformity: \\
There is a linear map 
$T: \hat{\scal}_o \rightarrow (\hat{\scal}_o)^n$, 
such that
\begin{enumerate}[\ \ \ (i)]
\item
 if $T\psi = (\psi_1,\ldots,\psi_n) := (T_1\psi,\ldots,T_n\psi)$, then 
   $\psi(\xi) = \sum_{j=1}^n \xi_j \psi_j(\xi)$; and
\item
 for every $I$ there exist $C, J$ such that for $1 \leq j \leq n$,
$\|T_j \psi\|_{\scal,I} \leq C\|\psi\|_{\scal,J}$.\\
    \end{enumerate}
\end{remark}
\begin{remark}
\label{rmk1L}
Taking inverse Fourier transforms in 
   Remark \ref{rmk1}, 
                for every $f \in \scal_o(G)$
   we can find $F_1,\ldots,F_n$ in $\scal_o$ 
   so that $f(x) = \sum_{j=1}^n \partial_j F_j(\xi)$.  
By Proposition \ref{proposicion-pagina-41}(b)
we can find $f_1,\ldots,f_n$ in $\scal_o$,
   with $f(x) = \sum_{j=1}^n X_j f_j(\xi)$. 
In fact, this construction shows the following uniformity: 
there exists a linear map $S: {\scal}_o \rightarrow ({\scal}_o)^n$, 
such that
\begin{enumerate}[\ \ \ (i)]
\item
    if $Sf = (f_1,\ldots,f_n) := (S_1f,\ldots,S_nf)$, then 
    $f(x) = \sum_{j=1}^n X_j f_j(x)$; and
\item
     for every $I$ there exist $C, J$ such that for $1 \leq j \leq n$,
    $$\|S_j f\|_{\scal,I} \leq C\|f\|_{\scal,J}.$$
\end{enumerate}
\end{remark}

\begin{remark}
\label{rmk2}
We can iterate the result of Remark \ref{rmk1}, by applying 
this result to each $\psi_j$ in place of $\psi$,
and so on.  We then obtain the following uniformity.  
For each $N \geq 1$, let $i(N)$ denote the number 
of multiindices $\alpha$ with $\|\alpha\| = N$.  
(Recall that the norm of a multiindex is equal to the sum 
of its coordinates, $\norm{\alpha}= \sum_{i=1}^n \alpha_i $. )
Then there is a linear map 
   $T^{(N)}: \hat{\scal}_o \rightarrow (\hat{\scal}_o)^{i(N)}$, such 
   that
\begin{enumerate}[\ \  \ (i)]
\item
     if $T^{(N)}\psi = (\psi_{\alpha})_{||\alpha||=N} 
          := (T^{(N)}_{\alpha}\psi)_{||\alpha||=N}$, then 
    $\psi(\xi) = \sum_{||\alpha||= N} \xi^{\alpha} \psi_{\alpha}(\xi)$; and
\item
     for every $I$ there exist $C$ and $J$ 
   such that for all $\alpha$ with $||\alpha||=N$, we have
   $$\relnorm{ T^{(N)}_{\alpha} \psi}_{\scal,I} \leq C\|\psi\|_{\scal,J}.$$
\end{enumerate}
\end{remark}

\begin{remark}
\label{rmk2L}
We can iterate the result of Remark \ref{rmk1L}, by applying this result 
  to each $f_j$ in place of $f$, and so on.  
We then obtain the following uniformity.  
For each $N \geq 1$, let ${\mathcal I}_N$ denote the collection 
  of $N$-tuples $\beta = (\beta_1,\ldots,\beta_N)$ 
with $1 \leq \beta_j \leq n$ for all $j$.  
For $\beta \in {\mathcal I}_N$, 
 let $X_{\beta} = X_{\beta_1}\cdots X_{\beta_{_N}}$.  
Note $\# {\mathcal I}_N = n^N$.  
Then there is a linear map 
   $S^{(N)}: \hat{\scal}_o \rightarrow (\hat{\scal}_o)^{n^N}$, such that
\begin{enumerate}[\ \  \ (i)]
\item
 if $S^{(N)}f = (f_{\beta})_{\beta \in {\mathcal I}_N} 
    := (S^{(N)}_{\beta}f)_{\beta \in {\mathcal I}_N}$, then 
    $f(x) = \sum_{\beta \in {\mathcal I}_N} X_{\beta}f_{\beta}(x)$; and
\item
 for every $I$ there exist $C$ and $ J$ such that for all 
   $\beta \in {\mathcal I}_N$, we have
  $$ \relnorm{ S^{(N)}_{\beta} f }_{\scal,I} \leq C\|f\|_{\scal,J}.$$
\end{enumerate}
\end{remark}

\section{\sc Characterization of $\mcalcop^j(G)$ } 
\begin{propo} \label{primera}
   If $m \in \mcal^j(G) $,
   there exists a bounded sequence
   $\relbrace{ \psi_k }_{k=0}^\infty \subset \scal(G)$, 
   with $\psi_k \in \hat{\scal}_o$
   for $ k \geq 1$, such that, if we set
  \begin{align} 
      m_k          &= 2^{jk} \psi_k \circ \delta_{2^{-k}},
     \label{ecuacion.1.primera}
      \\ 
\mbox{so that}\quad  
 \check{m}_k &= 2^{(j+Q)k }  \check{\psi_k} \circ \delta_{2^k} ,
     \label{ecuacion.2.primera}
  \end{align}
then $m = \sum_{k=0}^{\infty} m_k$, 
with convergence both pointwise and in $\scal'$.
\end{propo}

\begin{proof}
Choose a sequence of functions 
   $\{\varphi_k\}_{k=0}^{\infty} \subset C_c^{\infty}(G)$, such that:
\begin{enumerate}[\ \  \ (i)]
\item
 $\sum \varphi_k = 1$;
\item
 $0 \leq \varphi_k \leq 1$ for all $k$;
\item
 $\supp \varphi_o \subseteq \left\{ \xi : \abs{\xi} \leq 2 \right\}$;  
\item
 for all $k \geq 1$, $\supp \varphi_k \subseteq \left\{ \xi 
    : 2^{k-2} \leq \abs{\xi} \leq 2^{k+1}\right\}$; and
\item
 for all $k \geq 1$, $\varphi_k = \varphi_1 \circ \delta_{2^{-k+1}}$.
\end{enumerate}
(For example, we could choose a smooth function 
   $\varphi_o$ with $0 \leq \ \varphi_o\leq 1$, with $\varphi_o \equiv 1$
in $\left\{ \xi : \abs{\xi} \leq \frac{1}{2}\right\}$ and with 
$ \supp \varphi_o \subseteq  \left\{ \xi : \abs{\xi} \leq 2 \right\}$.  
  Then we could let $\varphi_1 = 
          \varphi_o \circ \delta_{1/2} - \varphi_o$, $\varphi_k 
    = \varphi_1 \circ \delta_{2^{-k+1}}$ for $k > 1$.)

Let $\varphi = \varphi_1 \circ \delta_2$, 
   so that $\varphi_k = \varphi \circ \delta_{2^{-k}}$
   for all $k \geq 1$, and so that $\supp \varphi \subseteq \left\{ \xi 
   : \frac{1}{4} \leq \abs{\xi} \leq 2\right\}$. 

We set $ m_k = \varphi_k m $  for $ k \geq 0 $, 
   so that $m = \sum_{k=0}^{\infty} m_k$, with pointwise
   convergence (at each point in $\rn$, 
   only finitely many terms are nonzero.)  (This is entirely analogous
to the definition used on page 246 of \cite{stein-libro-gordo}, for
standard multipliers $m$.)  We also define 
   $\psi_k$ so that (\ref{ecuacion.1.primera}) (and hence 
   (\ref{ecuacion.2.primera}))
   holds; thus $\psi_k = 2^{-jk} m_k \circ \delta_{2^k}$, and
   we therefore have
  \begin{align*} 
      \psi_o  & =   \varphi_o \left( m \circ  \delta_{2^k} \right) \\  
     \psi_k   & =   2^{-jk} \varphi \left(m \circ \delta_{2^k}  \right) 
                              \qquad \text{ for } k \geq 1.      
  \end{align*}
Since $\psi_k \in C_c^{\infty}$, and vanishes in a neighborhood of 
the origin if $k \geq 1$, surely $\check{\psi}_k \in \scal(G)$ 
for all $k$, and $\check{\psi}_k \in \scal_o(G)$ for $k \geq 1$.  

Note also that 
 $\supp \psi_k \subseteq \left\{ \xi : \abs{\xi} \leq 2 \right\}$ for all $k$.  
  Thus, to show that $\relbrace{ \psi_k }_{k=0}^\infty$ is a bounded subset 
  of $\scal(G)$, we need only show that for any multiindex 
   $\alpha$, $\relbrace{ \|\partial^\alpha \psi_k\|_{\infty} }_{k=1}^\infty$ 
  is a bounded sequence, where the sup norm is taken over 
   $\left\{ \xi : \abs{\xi} \leq 2 \right\}$.
But for $k \geq 1$, $\abs{\xi} \leq 2$, we have by Leibniz's rule that
  \begin{align*} 
      \relabs{ \partial^\alpha \psi_k(\xi)}   
         &\leq C_1 \sum_{ \abs{\beta} \leq \abs{\alpha} } 
                 2^{-jk} C_{\abs{\beta}} \ 2^{\abs{\beta}k} 
                \relabs{(\partial^{\beta} m) 
                    \relpar{ \delta_{2^k} (\xi) } } 
         \\
         &\leq C_2 \sum_{ \abs{\beta} \leq \abs{\alpha} } 
              2^{-( j -\abs{\beta} ) k } (1 + 2^k \abs{\xi} )^{j-\abs{\beta}}
         \\
         &= C_2 \sum_{ \abs{\beta} \leq \abs{\alpha} } 
          \left[  \frac{ 2^k }{(1+2^k\abs{\xi}) } \right]^{-(j-\abs{\beta})}
         \\
         & \leq C_\alpha
  \end{align*}
 as claimed.

Finally, note that, if $k \geq 1$, then for any $\xi \in \rn$,
$m_k(\xi) \neq  0$ for at most $3$ values of $k$. 
Moreover $|m_k(\xi)| = |\varphi_k(\xi)||m(\xi)| \leq C(1+|\xi|)^j$.  Thus
$\sum_k |m_k(\xi)| \leq C(1+|\xi|)^j$ also.   Since $\sum_k m_k$ converges
to $m$ pointwise, it now follows that the convergence is in $\scal'(G)$ 
as well.
\end{proof}

Next we show that the converse of the previous proposition also holds. 
We will also obtain more information about the senses in which 
$\sum_k m_k$ converges to $m$.
\begin{propo} \label{cnvchar}
\begin{enumerate}[ (a)]
\item
   Suppose  $\relbrace{ \psi_k }_{k=0}^\infty \subset \scal(G)$ 
   is a bounded sequence, with $\psi_k \in \hat{\scal}_o$ 
    for $k \geq 1$.  Set
  \begin{align*} 
      m_k  &   = 2^{jk} \psi_k \circ \delta_{2^{-k}}
      \\ 
    \text{so that}\phantom{mmm} \check{m}_k  &= 2^{(j+Q)k }  
    \check{\psi_k} \circ  \delta_{2^k}  ,
  \end{align*}
then $\sum_{k=0}^{\infty} m_k$ converges, 
   in $C^{\infty}(G)$, and in $\scal'$, to an element
$m$ of $\mcal^j(G)$.\\
\item  
 In fact, we have the following uniformity.  
   For every $I$ there exist $C$ and $J$ such that
whenever    $\relbrace{ \psi_k }_{k=0}^\infty \subset \scal(G)$ is 
   a bounded sequence, with
$\psi_k \in \hat{\scal}_o$ for $k \geq 1$, 
if we set $m_k = 2^{jk} \psi_k \circ \delta_{2^{-k}} $,
then
 $$ \relnorm{ \sum_{k=0}^{\infty} m_k}_{\mcal^j,I} 
     \leq C \sup_k \|\psi_k\|_{\scal,J}.  $$
\item
We also have the following uniformity.  Say $j' > j$.  
   Then, for every $I$ there exist $C$ and $J$ such that
   whenever    $\relbrace{ \psi_k }_{k=0}^\infty \subset \scal(G)$ 
   is a bounded sequence, with $\psi_k \in \hat{\scal}_o$ 
   for $k \geq 1$, if we set $m_k
   = 2^{jk} \psi_k \circ \delta_{2^{-k}}$,
then
 $$\relnorm{ \sum_{k=0}^{\infty} m_k  }_{\mcal^{j'},I} 
         \leq C \sup_k 2^{-(j'-j)k}\|\psi_k\|_{\scal,J}.  $$
\item
 In (a), the series $\sum_{k=0}^{\infty} m_k$ converges 
    to $m$ in $\mcal^{j'}(G)$, for any $j' > j$.
\end{enumerate}
\end{propo}
\begin{proof}
(a) Let us begin by showing that 
\begin{equation}\tag{$\ast$} \label{asterisk}
  \sum_k |m_k(\xi)| = \sum_k 2^{jk} 
   \relabs{  \psi_k (\delta_{2^{-k}} (\xi))} 
   \leq C(1+\abs{\xi})^j,\qquad
\end{equation}
with uniform convergence on compact sets; from this, 
 the convergence of $\sum_{k=0}^{\infty} m_k$, both pointwise 
 and in $\scal'$, will be automatic.  

First let us establish (\ref{asterisk}) for $\abs{\xi} > 1$.  
In that case, we define $k_0 \geq 1$ to be the unique integer 
  such that
   $2^{k_0 - 1} < \abs{\xi} \leq 2^{k_0}$.  
We consider separately those terms in 
   $\sum_k 2^{jk} \relabs{  \psi_k (\delta_{2^{-k}} (\xi))}$ 
with $k < k_0$ and with $k \geq k_0$.
Choose $N\in \z^+ $ such that $j+N > 0 $.
Since $\{\psi_k\}$ is a bounded subset of $\scal(G)$,
\begin{equation} \label{mayores.que.uno}
   \begin{split}
     \sum_{k < k_0} 
      2^{jk} \relabs{  \psi_k (\delta_{2^{-k}} (\xi)) } 
         & \leq C_1 \sum_{k < k_0 }
         2^{jk} (2^{-k} \abs{\xi} )^{-N}     			 \\
     & = C_1 \ \abs{\xi}^{-N} \sum_{k < k_0 }
           2^{ (j+N) k }  \\ 
     & \leq C_2 \abs{\xi}^{-N} 2^{ (j+N) k_0 }
   	\qquad  \text{ since }  j+N > 0 \\
	& \leq C_3 \abs{\xi}^j \\
          & \leq  \tilde{C_1} (1+\abs{\xi})^j
       \qquad\qquad \qquad  \qquad \text{ since }  \abs{\xi} > 1. 
   \end{split}
\end{equation}
Here $\tilde{C}_1$ may be chosed independently of $\xi$.  
 Note also that, for any fixed $k_0$, the convergence of
     $\sum_{k < k_0} 
      2^{jk} \relabs{  \psi_k (\delta_{2^{-k}} (\xi)) }$
is uniform for $2^{k_0 - 1} < \abs{\xi} \leq 2^{k_0}$.

Next, by Remark \ref{rmk2}, for every $k \geq 1 $ 
and for every $N \in \z^+$ we may write
 $\psi_k (\xi) 
  =\sum_{\|\alpha\|=N} \xi^\alpha \psi_{k,\alpha}(\xi),$ 
with $\psi_{k,\alpha}\in \scal(G)$, so that for any $N$, 
 $\{\psi_{k,\alpha}: k \geq 1, \|\alpha\|=N\}$
is a bounded subset of $\scal(G)$.  

  Now choose $N > j$.  Recall $k_0 \geq 1$.  We have
\begin{equation} \label{menores.que.uno}
   \begin{split}
     \hspace*{-5mm}
      \sum_{k \geq k_0} 
      2^{jk} \Big\vert
         \psi_k 
       \relpar{\delta_{2^{-k}} (\xi)} 
       \Big\vert
      & \leq
        C_1\sum_{k \geq k_0 } 
      2^{jk} \relpar{  
           \sum_{\|\alpha\| = N } 
              \relabs{ \delta_{2^{-k}} (\xi) }^{\abs{\alpha}}
                 \relabs{ \psi_{k,\alpha} \relpar{  \delta_{2^{-k}} (\xi)} } 
               } \\
      &\leq
       C_2 
        \sum_{k \geq k_0 }
          2^{jk} \relpar{  
              \sum_{\|\alpha\| = N } 
             {2^{-k{\abs{\alpha}}}} \abs{\xi }^{\abs{\alpha}}
                        }  \\
      &=
       C_3 \sum_{\|\alpha\| = N }
          \abs{\xi }^{\abs{\alpha}}
		\sum_{k \geq k_0 }
            2^{(j-{\abs{\alpha}})k} 
      \\ 
      &\leq 
          C_4 \sum_{\|\alpha\| = N }
          \abs{\xi }^{\abs{\alpha}} 2^{(j-{\abs{\alpha}})k_0}
      \\  &
       \qquad \text{ since } \|\alpha\| = N 
         \Rightarrow 
      j-{\abs{\alpha}} \leq j-N < 0\\
	& \leq C_5 {\abs{\xi}}^j\\
      & \leq
      \tilde{C}_2 (1+ \abs{\xi})^j
         \qquad \text{ since }  \abs{\xi} > 1.
   \end{split}
\end{equation}
Here again $\tilde{C}_2$ may be chosed independently of $\xi$.  Note also that, for any
fixed $k_0$, the convergence of
     $\sum_{k < k_0} 
      2^{jk} \relabs{  \psi_k (\delta_{2^{-k}} (\xi)) }$
is uniform for $2^{k_0 - 1} < \abs{\xi} \leq 2^{k_0}$.
Therefore by \eqref{mayores.que.uno} and 
     \eqref{menores.que.uno}, we have (\ref{asterisk}) for 
      $\abs{\xi} > 1$, with uniform convergence on bounded sets. 

Now if $\abs{\xi} \leq 1 $, using a similar argument as in 
 \eqref{menores.que.uno}, and again choosing $N> j$, we have
\begin{equation*}
\begin{split}
     \hspace*{-5mm}
      \sum_k 
      2^{jk} \Big\vert
         \psi_k 
       \relpar{\delta_{2^{-k}} (\xi)} 
       \Big\vert
	& =
         \relabs{ \psi_o (\xi) } 
        +
        \sum_{k=1}^\infty 2^{jk} \relabs{ \psi_k 
                                             (\delta_{2^{-k}} (\xi)) }
       \\
      & \leq
        C_1 + C_1\sum_k 
      2^{jk} \relpar{  
              \sum_{\|\alpha\| = N } \relabs{ \delta_{2^{-k}} (\xi) }^{\abs{\alpha}}
                 \relabs{ \psi_{k,\alpha} \relpar{  \delta_{2^{-k}} (\xi)} } 
               } \\
      &\leq
       C_1 + C_2
        \sum_{k=1}^\infty
          2^{jk} \relpar{  
              \sum_{\|\alpha\| = N } {2^{-k{\abs{\alpha}}}} \abs{\xi }^{\abs{\alpha}}
                        }  \\
      &=
       C_1 + C_2 \sum_{\|\alpha\| = N }
          		\sum_{k=1}^\infty
            2^{(j-{\abs{\alpha}})k} 
	\quad \text{ since } \abs{\xi}\leq 1 
      \\ 
      &\leq C
       \qquad \text{ since } \|\alpha\| = N \Rightarrow j-{\abs{\alpha}} \leq j-N < 0.\\
   \end{split}
\end{equation*}
Here $C$ may be chosen independently of $\xi$ (for $|\xi| \leq 1$).  Also the convergence
is uniform for $|\xi| \leq 1$.  This establishes (\ref{asterisk}), with uniform convergence on compact sets.

Finally, let $\alpha$ be a multiindex. Then 
\[    
\partial^\alpha m_k  = 2^{ (j-\abs{\alpha} )k } 
     \relpar{ \partial^\alpha \psi_k   } 
                  \circ \delta_{2^{-k}} . 
\]
Now, the sequence   
$\{ \partial^\alpha \psi_k \}_{k=0}^\infty$  is a bounded subset of $\scal(G)$,
  and $\partial^\alpha \psi_k \in \hat{\scal}_o$ for $k\geq 1$.
  Consequently, by (\ref{asterisk}), the series 
$\sum_k|\partial^\alpha m_k (\xi)|$ converges uniformly
on compact sets, and there exists a constant $C_\alpha$ such that 
\[ \sum_k|\partial^\alpha m_k (\xi)|  
       \leq
   C_\alpha (1+\abs{\xi} )^{j - \abs{\alpha} } \qquad \text{ for all } \xi. \]
Accordingly, the series $\sum_k m_k$ converges in $C^{\infty}(G)$, 
and in $\scal'(G)$, to an element $m \in \mcal^j(G)$.  This proves (a).

Part (b) follows from an examination of the proof of part (a) 
(and from the uniformity in Remark \ref{rmk2}).

Part (c) follows at once from (b), if we set 
   $\psi^{j'}_k = 2^{-(j'-j)k} \psi_k \circ\delta_{2^{-k}}$,
note that 
$\relbrace{ \psi^{j'}_k }_{k=0}^\infty \subset \scal(G)$ 
   is a bounded sequence, with
$\psi^{j'}_k \in \hat{\scal}_o$ for $k \geq 1$, and also note that
$m_k = 2^{j'k} \psi^{j'}_k \circ \delta_{2^{-k}}$.  

   Finally, part (d) follows at once from (c).
\end{proof}

We show now how the characterization of $\mcalcop^j (G)$
follows from the previous two propositions.
\begin{propo}
\label{mjcrt}
Say $K \in \scal'(G)$.  Then $K \in \check{\mcal}^j (G)$ 
  if and only if we may write
\begin{equation}
\label{kmgcrt}
K = \sum_{k=0}^\infty f_k,
\end{equation}
with convergence in $\scal'(G)$, where
\begin{equation}
\label{fkphi}
f_k = 2^{(j+Q)k} \varphi_k \circ \delta_{2^k}
\end{equation}
where    $\left\{\varphi_k \right\}_k \subseteq \scal(G)$  
 is a bounded sequence, and $\varphi_k \in \scal_o(G)$ for $k \geq 1$.    
Further, given such a sequence $\left\{\varphi_k \right\}_k$, 
if we define $f_k$ by (\ref{fkphi}), then the series $\sum_{k=0}^\infty f_k$ 
necessarily converges in $\scal'(G)$ to an element of $K$ 
of $\check{\mcal}^j (G)$.

Moreover, $\check{\mcal}^j (G) \subseteq {\mathcal E}' + \scal(G)$, 
and any $K \in \check{\mcal}^j (G)$ is smooth away from the origin.  
Further, if $K \in \check{\mcal}^j (G)$ is as in (\ref{kmgcrt}), and if 
$\zeta \in C_c^{\infty}(G)$ equals $1$ in a neighborhood of\, $0$, 
 then
\begin{equation}
\label{kmgcrts}
(1-\zeta)K = \sum_{k=0}^\infty (1-\zeta)f_k,
\end{equation}
where the sum converges absolutely in $\scal(G)$.
\end{propo}
\begin{proof}
The assertions of  the first paragraph are immediate 
from the previous two propositions. 
Thus, if $K \in \check{\mcal}^j (G)$ is as in (\ref{kmgcrt}), 
we have (\ref{kmgcrts}), with convergence in $\scal'(G)$; 
we need to show the sum converges absolutely in $\scal(G)$.  
For this, by Leibniz's rule, we only need show
that for every $I, J \geq 0$, 
         every multiindex $\alpha$,
         every $c > 0$, 
   and 
         every $r>0$, there exists $C > 0$ such that
$$ 
   \sum_{k=0}^{\infty} 2^{Ik}|x|^J 
   \relabs{  \partial^{\alpha} \varphi_k (\delta_{2^k}(x)) \phantom{ }
          } 
   \leq C 
$$
whenever $\abs{x} > r$. 

This, however, is evident, since for every $N, \alpha$ there 
exists $C_{N,\alpha}$ such that 
$$ |\partial^{\alpha} \varphi_k (\delta_{2^k}(x))| 
    \leq 
    C_{N,\alpha}(2^k|x|)^{-N}
    \qquad \text{ for all $x$ and all $k$.}
$$
The convergence of the sum in (\ref{kmgcrts}) in $\scal(G)$ 
  shows, in particular, that $(1-\zeta)K \in \scal(G)$, so
$K = \zeta K + (1-\zeta)K \in {\mathcal E}' + \scal(G)$.  
Finally, we may choose $\zeta$ to have support
within an arbitrarily small neighborhood of $0$.  
Outside this neighborhood, $K = (1-\zeta)K$ is smooth;
consequently, $K$ is smooth away from the origin.
\end{proof}
\section{\sc Convolving Elements of the ${\mcalcop}^j(G)$ } 
For $\sigma \in \z$, $I > 0$, define
\[ \Phi^I_{\sigma}(x) = (1+2^{\sigma}|x|)^{-I}. \]
We then have the following key fact.  
The statement and proof have been adapted 
from Lemma 3.3 of \cite{frazier-jawerth}
where only the case of 
Euclidean space is dealt with.
\begin{lemma}
\label{fjg}
Say $J > 0$, and let $I = J+Q$.
Then there exists $C > 0$ such that 
whenever $\sigma \geq \nu$, 
\begin{equation*} \label{fjgway}
    \Phi^I_{\sigma}*\Phi^J_{\nu} \leq C\, 2^{-\sigma Q} \Phi^J_{\nu}.
 \end{equation*}
\end{lemma}
\begin{proof}
We note that
\begin{align*} 
\Phi^I_{\sigma}*\Phi^J_{\nu}(x) 
& = \int_{|y| \leq |x|/2} \Phi^I_\sigma(y)\Phi^J_{\nu}(y^{-1}x) dy 
    + \int_{|y| \geq |x|/2} \Phi^I_\sigma(y)\Phi^J_{\nu}(y^{-1}x) dy \\
& := A + B. 
\end{align*}
In $A$ we note that, by (\ref{yxest}), 
   $\Phi^J_{\nu}(y^{-1}x) \leq C\, \Phi^J_{\nu}(x)$, so 
\[ 
   A \leq C\,\Phi^J_{\nu}(x)\int_G \Phi^I_\sigma(y)dy 
   = C\, 2^{-\sigma Q} \Phi^J_{\nu}(x), 
\]
since $I > Q$.  

In $B$, we just estimate $\Phi^J_{\nu}(y^{-1}x) \leq 1$.
Consider first the case $2^{\nu}|x| \leq 1$.  Then
\[ 
   B \leq \int_G \Phi^I_\sigma(y)dy 
   = C\,2^{-\sigma Q} \leq C\, 2^{-\sigma Q} \Phi^J_{\nu}(x). 
\]
Finally, if, instead, $2^{\nu}|x| \geq 1$, we have
\begin{equation*} 
  \begin{split}
  B \
  & \leq \ \int_{|y| \geq |x|/2} \Phi^I_\sigma(y)\ dy 
   \  \leq \ 
    \int_{|y| \geq |x|/2} (2^{\sigma}|y|)^{-I} dy \\
  & = C\,2^{-\sigma I}|x|^{-I+Q} 
   = C\,2^{-\sigma Q}(2^{\sigma}|x|)^{-J}. 
 \end{split}
\end{equation*} 
Thus, since $\sigma \geq \nu$,
\begin{equation*} 
   B 
   \leq 
      C\,2^{-\sigma Q}(2^{\sigma}|x|)^{-J} \\
  \leq C\, 2^{-\sigma Q}(2^{\nu}|x|)^{-J} \\
  \leq C\,2^{-\sigma Q} \Phi^J_{\nu}(x), 
\end{equation*} 
as desired.
\end{proof}

\begin{remark}
One may think of Lemma \ref{fjg} very crudely in the following manner: from the perspective of
$\Phi^J_{\nu}$, $2^{\sigma Q} \Phi^I_{\sigma}$ looks like a good approximation to the delta ``function'', so 
convolving with it gives one back something akin to $\Phi^J_{\nu}$.
\end{remark}

\begin{lemma}
\label{bdcnvbd}
Suppose $L \geq 0$ is an integer, and that ${\mathcal B} \subseteq \scal(G)$ 
 is bounded.  Then there is a bounded subset 
   ${\mathcal B}' \subseteq \scal(G)$ as follows.
Say $k, l \in \z$, $k \geq l$.  Suppose that 
   $\varphi \in {\scal}_o(G) \cap {\mathcal B}$, and that
$\psi \in {\mathcal B}$.  Define $w_{k,l} \in \scal_o(G)$ by 
\begin{equation} \label{hklfg}
  w_{k,l} \circ \delta_{2^l} 
   = 2^{kQ} (\varphi \circ \delta_{2^k}) * (\psi \circ \delta_{2^l}).
\end{equation}
Then $2^{(k-l)L}w_{k,l} \in {\mathcal B}'$.
\end{lemma}
\begin{proof}
This proof has three steps.\\
{\sc Step 1:\ } 
Say $J > 0$, and let $I = J+Q$.  Then by Lemma \ref{fjg},
\[ 
|w_{k,l} \circ \delta_{2^l}| \leq C\, 2^{kQ}\Phi^I_{k}* \Phi^J_{l} 
               \leq  C\, \Phi^J_{l} = C\, \Phi^J_0 \circ \delta_{2^l}, 
\]
so $|w_{k,l}| \leq C\, \Phi^J_0$.  In other words,
\[ 
   (1+|x|)^J w_{k,l}(x) \leq C. 
\]
{\sc Step 2:\ } 
Notation as in Remark \ref{rmk2L}, say $\beta \in {\mathcal I}_N$, 
and let $r = \sum_{j=1}^{_N} a_{\beta_j}$.  Applying $X_{\beta}$ 
to both sides of (\ref{hklfg}), we see that
\[ 
   2^{rl} (X_{_\beta}w_{k,l}) \circ \delta_{2^l} = 2^{kQ} (\varphi \circ \delta_{2^k}) * [2^{rl}(X_{\beta}\psi) \circ \delta_{2^l}],
\]
so that
\[ 
(X_{\beta}w_{k,l}) \circ \delta_{2^l} 
  = 2^{kQ} (\varphi \circ \delta_{2^k}) * [(X_{\beta}\psi) \circ \delta_{2^l}].
\]
Note that $X_{\beta}\psi \in X_{\beta}{\mathcal B}$, 
a bounded subset of $\scal(G)$.  Thus, by
Step 1, for any $J > 0$ and any $\beta$, we have
\[ (1+|x|)^J (X_{\beta}w_{k,l})(x) \leq C_{J,\beta}. \]
This proves the lemma in the case $L=0$.\\
{\sc Step 3:\ } Finally, suppose $L \geq 1$.  
By Remark \ref{rmk2L}, we may assume that for some bounded subset
${\mathcal B}_L$ of $\scal(G)$ (depending only on ${\mathcal B}$ and $L$), 
$\varphi = X_{\beta}\varphi_1$, where
$X_{\beta} \in {\mathcal I}_L$, and
$\varphi_1 \in {\scal}_o(G)\cap{\mathcal B}_L$.  
Set $Y_{\beta} = Y_{\beta_{_L}}\cdots Y_{\beta_{_1}}$.  
Also set $r = \sum_{j=1}^L a_{\beta_j} \geq L$.  Now 
\[ 
   \varphi \circ \delta_{2^k} 
   = (X_{\beta}\varphi_1) \circ \delta_{2^k} 
   = 2^{-rk} X_{\beta}(\varphi_1 \circ \delta_{2^k}), 
\]
so
\begin{equation*}
\begin{split}
w_{k,l} \circ \delta_{2^l} & = 2^{kQ} [ 2^{-rk} X_{\beta}(\varphi_1 \circ \delta_{2^k})] * (\psi \circ \delta_{2^l})\\
& = 2^{kQ} 2^{-rk} (\varphi_1 \circ \delta_{2^k}) * Y_{\beta}(\psi \circ \delta_{2^l})\\
& = 2^{-(k-l)r}2^{kQ}  (\varphi_1 \circ \delta_{2^k}) * [(Y_{\beta}\psi) \circ \delta_{2^l}].
\end{split}
\end{equation*}
Since $r \geq L$, we therefore have that
\[ (2^{(k-l)L} w_{k,l}) \circ \delta_{2^l} 
   = c\ 2^{kQ}  (\varphi_1 \circ \delta_{2^k}) * 
       [(Y_{\beta}\psi) \circ \delta_{2^l}], 
\]
where $0 < c \leq 1$.  Note also that $Y_{\beta} \psi \in Y_{\beta}{\mathcal B}$, a bounded subset
of $\scal(G)$.  By the case $L=0$ of the lemma, established in Step 2, we see that there is a bounded
subset ${\mathcal B}'$ of $\scal(G)$, depending only on ${\mathcal B}$ and $L$,
such that $2^{(k-l)L}w_{k,l} \in {\mathcal B}'$. 
\end{proof}

We now reformulate Lemma \ref{bdcnvbd} in a manner that permits us to deal with the 
convolution of sums like that in (\ref{ecuacion.2.primera}).

\begin{coro}
\label{cnvfkgk}
Suppose $L \geq 0$ is an integer, and that ${\mathcal B} \subseteq \scal(G)$ is bounded.  Then there is a bounded
subset ${\mathcal B}' \subseteq \scal(G)$ as follows.
Say $k, l \in \z$, $k \geq l$, and that $j_1, j_2$ are real numbers.  
Suppose that $\varphi \in {\scal}_o(G) \cap {\mathcal B}$, and that
$\psi \in {\mathcal B}$.  Define $f \in \scal_o(G)$, $g \in \scal(G)$, $\eta_{k,l}, \eta^{\prime}_{k,l} \in \scal_o(G)$ by 
  \begin{align*} 
  f  & = 2^{(j_1 +Q ) k}\  \varphi \circ \delta_{2^{k}}  \\
  g  & = 2^{(j_2 +Q )l} \ \psi \circ \delta_{2^{l}}\\ 
  f*g& = 2^{(j_1+j_2+Q)l}\  \eta_{k,l} \circ \delta_{2^l}\\
  g*f& = 2^{(j_1+j_2+Q)l}\  \eta^{\prime}_{k,l} \circ \delta_{2^l}	
  \end{align*}
Then $2^{(k-l)(L-j_1)}\ \eta_{k,l} \in {\mathcal B}'$, and 
   also $2^{(k-l)(L-j_1)}\ \eta^{\prime}_{k,l} \in {\mathcal B}'$. 
  \end{coro}
\begin{proof}
We have $2^{-(k-l)j_1}\eta_{k,l} = w_{k,l}$, where $w_{k,l}$ is as in (\ref{hklfg}).  
Thus the first statement follows at once from Lemma \ref{bdcnvbd}.  For the second statement, we need only
apply the first statement to $\tilde{f}, \tilde{g}$ in place of $f,g$, and then use the fact that
$g*f = (\tilde{f}*\tilde{g})\tilde{\ \ }$. 
\end{proof}

We are almost ready to prove our main result, that 
    $
    \check{\mcal}^{j_1}
     \ast
    \check{\mcal}^{j_2}
    \subseteq
    \check{\mcal}^{j_1+j_2}
    $.
Before we begin the proof, it will be helpful to make 
some preliminary observations.

Say $K_1 \in \check{\mcal}^{j_1}(G)$, $K_2 \in \check{\mcal}^{j_2}(G)$.  
By Proposition \ref{mjcrt}, $K_1, K_2 \in {\mathcal E}' + \scal(G)$, 
which is a {\em convolution algebra}; thus it is possible to form 
$K_1 \ast K_2$.  (Of course we cannot convolve two general elements of 
$\scal'(G)$.)  To see that ${\mathcal E}' + \scal(G)$ is a convolution algebra,
 say $F, G \in {\mathcal E}' + \scal(G)$.
To convolve them we simply write $F = u + f$, $H = v + g$, 
where $u,v \in {\mathcal E}'$ and $f,g \in \scal(G)$; and then we define
\begin{equation}
\label{fastgdf}
F \ast H = u \ast v + u \ast g + f \ast v + f \ast g.
\end{equation}
The first term on the right side is in ${\mathcal E}'$ and the 
other terms are in $\scal(G)$, 
  so $F \ast H$ is in ${\mathcal E}' + \scal(G)$.  
It is easy to verify that this definition of $F \ast H$ is independent
of how one chooses to decompose $F$ in the form $u + f$, 
or $H$ in the form $v + g$.

It is also evident from (\ref{fastgdf}) that if $u_N \rightarrow u$ and $v_N \rightarrow v$ in ${\mathcal E}'$,
while $f_N \rightarrow f$ and $g_N \rightarrow g$ in $\scal(G)$, then
\begin{equation}
\label{umfm}
(u_N + f_N) \ast (v_N + g_N) \longrightarrow (u + f) \ast (v + g)
\end{equation}
in $\scal'(G)$.

We can now prove our main theorem:
\begin{teo} 
 If $j_1,j_2$ are real numbers then 
    $
    \check{\mcal}^{j_1}(G)
     \ast
    \check{\mcal}^{j_2}(G)
    \subseteq
    \check{\mcal}^{j_1+j_2}(G)
    $
\end{teo} 
\begin{proof}
Say $K_1 \in \check{\mcal}^{j_1}(G)$, and $K_2 \in \check{\mcal}^{j_2}(G)$.  
As in Proposition \ref{primera}, we may write 
\[ 
    K_1 =  \sum_{k_1=0}^\infty f_{k_1},\qquad
    K_2 =  \sum_{k_2=0}^\infty g_{k_2}, \]
(convergence in $\scal'(G)$), where
\begin{align*}
f_{k_1} &= 2^{ (j_1 + Q)k_1} \varphi_{k_1} \circ \delta_{2^{k_1}} 
& g_{k_2} &= 2^{ (j_2 + Q)k_2} \psi_{k_2} \circ \delta_{2^{k_2}};
\end{align*} 
where $\{ \varphi_{k} \} $ and $ \{\psi_{k}\}$ 
are bounded sequences in $\scal(G)$, and
$\varphi_k, \psi_k \in \scal_o (G) $ for $k\geq 1$.

We want to write $K_1 * K_2$ in the form  
$\sum_k 2^{(j_1+j_2+Q)k }\nu_k \circ \delta_{2^k}$ (convergence in $\scal'(G)$),
where $\{ \nu_{k} \}$ is a bounded sequence in $\scal(G)$, and where
$\nu_k \in \scal_o (G) $ for $k\geq 1$.  
Then we will know $K_1*K_2 \in \check{\mcal}^{j_1+j_2}(G) $ by
Proposition \ref{cnvchar}.

To this end we examine $f_{k_1}*g_{k_2} := h_{k_1,k_2}$.   
Formally,  $K_1*K_2 = \sum_{k_1,k_2} h_{k_1,k_2}$; 
we intend to write the double summation here as a sum
of two double summations, one for $k_2 \geq k_1$, and one for $k_1 > k_2$.  
Note that, to say that
$k_2 \geq k_1$, is to say that either 
$k_2 \geq \max(k_1,1)$, or $(k_1,k_2) = (0,0)$.

Pick $L > \max(j_1,j_2)$.  By Corollary \ref{cnvfkgk}, 
   we have 
\begin{align*}
h_{k_1,k_2} 
    &= 2^{(j_1+j_2+Q)k_1} \tau_{k_1,k_2} \circ \delta_{2^{k_1}}
    &&  \text{ for } k_2 \geq \max(k_1,1)\\
 h_{k_1,k_2} &= 2^{(j_1+j_2+Q)k_2} \eta_{k_1,k_2} \circ \delta_{2^{k_2}}
   && \text{ for } k_1 > k_2.
\end{align*}
where 
$ \left\{2^{(k_2-k_1)(L-j_2)} \tau_{k_1,k_2}\right\}_{k_2 \geq \max(k_1,1)}$ 
  and
$\left\{2^{(k_1-k_2)(L-j_1)}\eta_{k_1,k_2}\right\}_{k_1 > k_2}$ 
are bounded subsets of $\scal_o(G)$.

Formally, then,
\begin{eqnarray}
   K_1*K_2 &= & h_{o,o} + 
           \hspace*{-8mm}
      \sum_{\hspace*{10mm} k_2 \geq \max(k_1,1)} 
                    \hspace*{-4mm}
               h_{k_1,k_2} 
              + \sum_{k_1 > k_2} h_{k_1,k_2} 
                \label{k1k21}
             \\
           &= & h_{o,o} 
               + \sum_{k_1 = 0}^{\infty}
                            \hspace*{-5mm}
                            \sum_{\hspace*{7mm} k_2 \geq \max(k_1,1)}^\infty
                             \hspace*{-4mm}
                 h_{k_1,k_2} +
                   \sum_{k_2 = 0}^{\infty}
                   \hspace*{-5mm}
                   \sum_{\hspace*{7mm} k_1 = k_2 + 1}^{\infty} 
                   \hspace*{-4mm} h_{k_1,k_2} 
                    \label{k1k22}
                    \\
           &= & h_{o,o} + \sum_{k_1=0}^{\infty}u_{k_1} 
                + \sum_{k_2=0}^{\infty}v_{k_2} \label{k1k23}
\end{eqnarray}
where
\begin{align}
   u_{k_1} &= 2^{(j_1+j_2+Q)k_1} \tau_{k_1} \circ \delta_{2^{k_1}} \notag
  & v_{k_2} &= 2^{(j_1+j_2+Q)k_2} \eta_{k_2} \circ \delta_{2^{k_2}}  \\
\intertext{ and where }
  \tau_{k_1} &= \sum_{k_2 = \max(k_1,1)}^{\infty}\tau_{k_1,k_2} 
  & \eta_{k_2} &= \sum_{k_1 = k_2+1}^{\infty}\eta_{k_1,k_2}. \label{tauk}
\end{align}
The point is that these series  converge absolutely in 
$\scal(G)$, and that $\{\tau_{k_1}\}$, $\{\eta_{k_2}\}$ 
are bounded sequences in $\scal_o(G)$.  
(Once this is verified, we will know, at least formally, that 
$K_1*K_2 \in \check{\mcal}^{j_1+j_2}$.)  
But this point is easily verified.  
 Let $\norm{\cdot}_{\scal,M}$ be 
  any member of the family of norms defining the topology of $\scal(G)$.
Then, for instance
in the first series of~(\ref{tauk}) we have
\begin{align} 
\sum_{k_2 = \max(k_1,1)}^{\infty}\|\tau_{k_1,k_2}\|_{\scal,M} 
  & \leq C_M 
      \hspace*{-9mm}
    \sum_{\hspace*{8mm} k_2 = \max(k_1,1)}^\infty
   \hspace*{-4mm}
   2^{-(k_2-k_1)(L-j_2)} \notag  \\
  & \leq C_M\sum_{k=0}^{\infty} 2^{-k(L-j_2)} \notag \\
  &  = C_M^{\prime}\label{tauabcv}
\end{align}
where $C_M$ and $C_M^{\prime}$ are {\em independent} of $k_1$.  
Similar considerations apply to the the other series.
Of course the fact that $\tau_{k_1}, \eta_{k_2} \in \scal_o(G)$ 
is ensured by the absolute convergence of the series in $\scal(G)$ 
in (\ref{tauk}). 

Therefore we only need verify that $K_1*K_2$ is in fact equal 
to (\ref{k1k23}), as elements of $\scal'(G)$.  By Proposition 
\ref{mjcrt}, we do know that (\ref{k1k23}) converges in 
$\scal'(G)$.

Let us first show that
\begin{equation}
\label{k1k2lim}
\left(\sum_{k_1=0}^{_N} f_{k_1} \right) 
   * \left(\sum_{k_2=0}^{_N} g_{k_2} \right) \longrightarrow K_1*K_2
\end{equation}
in $\scal'(G)$, as $N \rightarrow \infty$.
Indeed, choose $\zeta \in C_c^{\infty}(G)$ with $\zeta = 1$ near $0$.  
We may write
\[ \sum_{k_1=0}^{_N} f_{k_1} 
  = \zeta \left(\sum_{k_1=0}^{_N} f_{k_1} \right) 
     + (1 - \zeta) \left(\sum_{k_1=0}^{_N} f_{k_1} \right). \]
As $N \rightarrow \infty$, we have that
   $\zeta \left(\sum_{k_1=0}^{_N} f_{k_1}\right) \rightarrow \zeta K_1$ 
   in ${\mathcal E}'$,
while, by Proposition \ref{mjcrt}, 
 $$(1 - \zeta) \left(\sum_{k_1=0}^{_N} f_{k_1}\right) 
                  \rightarrow (1-\zeta)K_1$$
in $\scal(G)$.  Similar considerations apply to 
$\sum_{k_2=0}^{_N} g_{k_2}$.  Thus (\ref{k1k2lim})
follows from (\ref{umfm}).

Thus we need only show that 
  $\left(\sum_{k_1=0}^{_N} f_{k_1}\right) * 
   \left(\sum_{k_2=0}^{_N} g_{k_2}\right)$ approaches
(\ref{k1k23}) in ${\mathcal S}'$, as $N \rightarrow \infty$.  

Let us make the convention that any sum of the form $\sum_{k=i}^j$ 
is zero if $j < i$.
In place of (\ref{k1k22}) and (\ref{k1k23}), we then rigorously have that 
\begin{align*}
\left(\sum_{k_1=0}^{_N} f_{k_1}\right) 
        * \left(\sum_{k_2=0}^{_N} g_{k_2}\right)
 &= h_{o,o} + \sum_{k_1 = 0}^{_N} 
   \hspace*{-3mm}
    \sum_{\hspace*{7mm} k_2 = \max(k_1,1)}^{_N} 
   \hspace*{-3mm}
             h_{k_1,k_2} 
   + \sum_{k_2 = 0}^{_N}
      \sum_{\hspace*{3mm} k_1 = k_2+1}^{_N} h_{k_1,k_2} \\
 &= h_{o,o} + \sum_{k_1=0}^{_N}u^N_{k_1}+ \sum_{k_2=0}^{_N}v^N_{k_2},
\end{align*}
where
\begin{align*}
u^N_{k_1} &= 2^{(j_1+j_2+Q)k_1} \tau^N_{k_1} \circ \delta_{2^{k_1}}
& v^N_{k_2} &= 2^{(j_1+j_2+Q)k_2} \eta^N_{k_2} \circ \delta_{2^{k_2}}
\intertext{and where }
\tau^N_{k_1} &= \sum_{k_2 = \max(k_1,1)}^{N}\tau_{k_1,k_2}
& \eta^N_{k_2} &= \sum_{k_1 = k_2+1}^{N}\eta_{k_1,k_2}.
\end{align*}

The argument of (\ref{tauabcv}) shows 
that $\{\tau^N_{k_1}\}_{k_1,N}$ and $\{\eta^N_{k_2}\}_{k_2,N}$
are bounded subsets of $\scal_o(G)$, 
and moreover that $\tau^N_{k_1} \rightarrow \tau_{k_1}$ in
$\scal(G)$ for each $k_1$, 
and $\eta^N_{k_2} \rightarrow \eta_{k_2}$ in $\scal(G)$ 
for each $k_2$, as $N \rightarrow \infty$.

To complete the proof, we need to show that 
\[ 
\sum_{k_1=0}^{N}u^N_{k_1} \longrightarrow \sum_{k_1=0}^{\infty}u_{k_1} 
  \qquad \text{ in } \scal'(G)  
\]
and that
\[ 
\sum_{k_2=0}^{N}v^N_{k_2} \longrightarrow \sum_{k_2=0}^{\infty}v_{k_2} 
  \qquad \text{ in } \scal'(G)  .
\]
Let us prove the first of these; the proof of the second 
is virtually identical.

We work on the Fourier transform side, and we write $k$ 
in place of $k_1$ for simplicity.  
Set $m_k = \hat{u}_k$ and $m^N_k = \hat{u}^N_k$.  
  We need to show that
\[ 
\sum_{k=0}^{N} m^N_{k} \longrightarrow \sum_{k=0}^{\infty}m_{k} 
  \qquad \text{ in } \scal'(G). 
\]
Set $\psi_k = \hat{\tau}_{k}$ and $\psi^N_k = \hat{\tau}^N_{k}$.  
  We have that
\[ m^N_k = 2^{(j_1+j_2)k} \psi^N_{k} \circ \delta_{2^{-k}}\]
\[ m_k = 2^{(j_1+j_2)k} \psi_{k} \circ \delta_{2^{-k}}\]
that $\{\psi^N_{k}\}_{k,N}$ and $\{\psi_{k}\}_{k}$
are bounded subsets of $\hat{\scal}_o$, 
and moreover that $\psi^N_{k} \rightarrow \psi_{k}$ 
in $\scal(G)$ for each $k$,  as $N \rightarrow \infty$.  

Proposition \ref{cnvchar} (c) now implies at once 
that $\sum_{k=0}^{N} m^N_{k} \rightarrow \sum_{k=0}^{\infty}m_{k}$
in $\mcal^{j'}(G)$ for any $j' > j_1+j_2$.  In particular, 
$\sum_{k=0}^{N} m^N_{k} \rightarrow \sum_{k=0}^{\infty}m_{k}$ 
in $\scal'(G)$, as desired.
\end{proof}
\section{\sc Pseudodifferential Calculus on $G$}
\subsection{\sc introduction}
\begin{defi}
 We say that an operator $A$ in $\scal(\rn)$ 
 is a {\em Fourier multiplier operator} if 
for a suitable function $a$
\begin{equation}\label{seis.uno.primera} 
 \big[ A f\big] (x) 
 = \int_\rn e^{-2\pi i x\cdot \xi} a(\xi)  \, \widehat{f}(\xi) \, d\xi
    \end{equation}
for $f\in\scal(\rn)$.
The function $a(\xi)$ is called the {\em multiplier}. 
\end{defi}

 Notice that Fourier multiplier operators can be written
 as convolution operators
  \begin{equation*}\label{} 
  \big[  A f   \big] (x) = \big(\,a \widehat{f} \,\big)\spcheck(x) 
                = \big( \copete{a}\ast f)(x)
    \end{equation*}
 where $\ast$ is here the usual convolution on $\rn$.
 $\sro{\#}$ will be used to denote the space of those Fourier multipliers
  which are in $\sro{0}$. For $p\in\sro{\#}$ 
   we define a family of seminorms by
  \begin{equation*}\label{} 
    \norm{p}_{\alpha,m,\rho} = \sup_{\xi} 
           \left\{ \gul{\xi}^{-m+\rho \norm{\alpha}  }
            \bigabs{ D_\xi^\alpha p(\xi)  }
           \right\}
    \end{equation*}
 where $\gul{\xi}= (1+\norm{\xi}^2)^{1/2}$. Equipped with these seminorms
 $\sro{\#}$ is a Fr\'echet space.

Recall that given $J\in\mcal^j(\rn)$ 
the \emph{multiplier operator} $ {m_J} $ is defined by 
         \begin{equation*}\label{}
    \begin{split}
   \hspace*{4mm}\xymatrix@1{
        \scal' \ \ar@{->}[rr]^{m_J}
           & & {\ \hspace{1mm} \scal'\hspace{3mm}}
              } \\
       \xymatrix@1{
    {\hspace{0mm}   \ f   \ }
                  \ar@{|->}[rr]
           &   &{\ \copete{J}  \ast f }
              }
      \end{split}
    \end{equation*}
Note $\mcal^j(\rn) = S_{1,\#}^j$.

Suppose  $a(x,D)$   is in $ \Op{\sro{0}}$ 
with  symbol   $a(x,\xi)$. For any $x\in \rn$ we define 
$a_x(\xi) = a(x,\xi)$.
If  $A(y)$ is the operator of Fourier multiplication by
 $a_y(\xi)$, then for any $f\in \scal(\rn)$
\begin{equation}\label{cinco-dos} 
    \begin{split}
  \big[ a(x,D) f\big] (x) 
    & = \int e^{-2\pi i x\cdot \xi}\ a(x,\xi) \ \widehat{f}(\xi) \ d\xi \\
    & = \int e^{-2\pi i x\cdot \xi}\ a_x(\xi)  \ \widehat{f}(\xi) \ d\xi \\     
    &= \big[ A(x) f\big] (x) 
      \end{split}
  \end{equation}
Therefore, loosely speaking, one can say that
a pseudodifferential
 operator is a {\em multiplier operator where the multiplier
 depends smoothly upon the point at which we are}.

Since 
$ [ A(x)f] (x) = (\copete{a}_x \ast f ) (x)$
then 
\begin{equation*}\label{} 
[ a(x,D)f] (x) = (\copete{a}_x \ast f ) (x) .
    \end{equation*}
Therefore, locally one can always represent 
 a pseudodifferential operator on $\rn$ with symbol in $S_{\rho,0}$, 
  by a smooth family of convolution operators, 
 where one 
 convolves with an element $\copete{a}_x \in \big(\sro{\#}\big)\spcheck$.
 This point of view is useful when working with pseudodifferential operators 
 on homogeneous groups $G$, where the group Fourier transform is cumbersome
to use.
Taylor shows in \cite{taylor-part-i} that 
smooth families of convolution operators can also be used 
to construct certain classes of pseudodifferential operators 
on $G$.  
 If $G$ is a homogeneous group, instead of Fourier multiplier operators
 as in (\ref{seis.uno.primera}), one 
 considers convolution operators on $G$ , defined by
  \[      A f = 
        \copete{a} \ast f  \]
     where $\ast$ is now group convolution. One requires that
 $a$ belong to some Fr\'echet space 
 ${\frakX}$ 
  of smooth functions on $\rn$.
 One assumes that $\frakX \subseteq \sro{\#}$ 
 for some $ m \in \rr $ and $\rho \in (0,1]$, and says that 
  $A\in \Op{\frakX}$.

      Say now that, instead, $a(y,\xi) \in \ci (U\times\rn)$
      where $U\subseteq G$ is open.
      Set $a_y (\xi) = a(y,\xi)$ and suppose
      $a_y$ is a smooth function of $y$, taking values in
      $\frakX$, for $y\in U$.
      For $y \in U$, one defines an operator 
      $A(y): \ci_c(G) \longrightarrow \ci(G) $
     by 
\begin{equation} \label{ayfchkdf}
      A(y) f  =  \check{a}_y \ast f
\end{equation}
     Then for $y\in U$ one defines 
\begin{equation} \label{frkafdf}
     \left[ \frakA f \right] (y) =  \left[ A(y) f \right] (y) 
\end{equation}
  Notice the analogy with (\ref{cinco-dos}).
  One denotes the set of such operators by $\Op{\widetilde{\frakX}}$.
  When $a(y,z)$ has compact support in $y$, we denote the collection of such
  operators $\frakA$ by $\Opc{\widetilde{\frakX}}$.

The following theorem contains results from  Taylor \cite{taylor-part-i}.
\begin{teo}\label{teorema-taylor}
   Suppose $G$ is a Lie group, and 
   $\displaystyle \left\{ \frakX^m \right\}_{m\in\rr}$
   is a nested family of Fr\'echet spaces 
   satisfying the following  conditions
 \begin{enumerate}[\ \ \ \ \  a) \  ]
  \item  \label{hypo-uno}  
        If $m\geq 0$ then    
        $\displaystyle \frakX^m \subset \sro{\#}$ 
         for some $\rho \in (0,1]$. 
  \item  If $m< 0$ then    
        $\displaystyle \frakX^m \subset \sroo{m\sigma}{\#}$ 
         for some $\sigma \in (0,1]$. 
  \item 
     If 
     $A\in \Op{\frakX^{m_1}}$, and 
     $B\in \Op{\frakX^{m_2}}$, 
     then $AB \in  \Op{\frakX^{m_1+m_2}}$,
     the product being continuous.
  \item    \label{hypo-cuatro}
    If $p(\xi)\in \frakX^m$, then 
    $D_\xi^\alpha p(\xi)\in \frakX^{m-\tau\|{\alpha}\|}$ for
    some $\tau \in (0,1]$.
  \item 
   If $K_j\in \frakX^{m-\tau j}$, then there exists 
    $K\in \frakX^m $ such that, for any $M$,
     if $N$ is sufficiently large, 
    \[   K - \left( K_o + \cdots + K_N\right) \in S^\peq{-M}_\peq{\rho,\#}\]
  \item 
   If $p(\xi) \in \frakX^m$ then $\bar{p}(\xi) \in \frakX^m$.
 \end{enumerate}
Then on $\Op{\widetilde{\frakX}}$ we have a pseudodifferential 
calculus, of the {\em usual type}
 containing products and adjoints, more specifically  
\begin{enumerate}[\ \ i)]
\item \label{thesis-i} 
      If\, $\frakA\in\Op{\widetilde{\frakX}^{m_1}}$ 
      (as in (\ref{ayfchkdf}), (\ref{frkafdf})),  
       and\ $\frakB\in\Opc{\widetilde{\frakX}^{m_2}}$,
        \\
      then $\frakA\frakB \in \Op{\widetilde{\frakX}^{m_1 + m_2}}$, 
      and it has an asymptotic expansion
  \begin{equation} \label{ABasymp}
      \big [\frakA \frakB f \big] (x)
      \sim
      \sum_{\gamma\geq 0}  
           \big [A^{\peq{[\gamma]}} (x) 
           B_\peq{{[\gamma]}}(x)f \big](x),
  \end{equation}
   in the sense that the operator 
   $\frakA \frakB 
   -  \sum_{\|\gamma\| \leq N} A^{\peq{[\gamma]}} (x) B_\peq{{[\gamma]}}(x)$
    becomes arbitrarily highly smoothing as $N \rightarrow \infty$.  
   Here the operators  $A^{\peq{[\gamma]}} (x)$,
   $B_\peq{{[\gamma]}}(x)$ are of the form
 \begin{equation*} \label{AgamBgam}
    A^{\peq{[\gamma]}} (x)g = \check{a}^{\peq{[\gamma]}}_x*g
        \qquad  \qquad  
    B_{\peq{[\gamma]}} (x)g = \check{b}_{\peq{[\gamma]},x}*g
 \end{equation*}
    where
   \begin{equation*} \label{agamdf}
   a^{\peq{[\gamma]}}_x(\xi) = (D_{\xi}^{\gamma}a)(x,\xi)
   \end{equation*}
and
\begin{equation*} \label{bgamprp}
   b_{\peq{[\gamma]},x}(\xi) 
   := 
   b_{\peq{[\gamma]}}(x,\xi) \in \ci(\rn \times \rn) 
           \mbox{ and is compactly supported in } x,
\end{equation*}
 and moreover
\begin{equation*} \label{bgamdf}
   b_{\peq{[\gamma]},x} 
        \mbox{ is a smooth function of } x, 
        \mbox{ with values in } \frakX^{m_2}.   
\end{equation*}
Then by hypothesis \ref{hypo-cuatro})  of this theorem, 
    $A^{\peq{[\gamma]}} (x) B_\peq{{[\gamma]}}(x) 
     \in 
     \Op{\widetilde{\frakX}^{m_1 + m_2 - \tau \|\gamma\|}}$, 
  so that (\ref{ABasymp}) is an asymptotic expansion within the 
  $\Op{\widetilde{\frakX}^{\mu}}$ spaces.
\item   \label{thesis-ii} 
    If $\frakA\in\Op{\widetilde{\frakX}^{m}}$ then the adjoint 
     $\frakA^*\in \Op{\widetilde{\frakX}^{m}}$ and it has 
     the following asymptotic expansion
\begin{equation} \label{Aastasy}
    \big[ \frakA^*  f \big] (x) 
    \sim
    \sum_{{\gamma}\geq 0}  
          \big [A^{\peq{\{\gamma\}}} (x) 
               f \big](x)
    \end{equation}
   in the same sense as in \ref{thesis-i}).
   Here the operators $A^{\peq{\{\gamma\}}} (x)$ are of the form
\begin{equation*} \label{Aastgam}
    A^{\peq{\{\gamma\}}} (x)g 
      = \check{a}^{\peq{\{\gamma\}}}_x*g
\end{equation*}
   where
   \begin{equation*} \label{aastc}
        a^{\peq{\{\gamma\}}}_x(\xi) 
        = 
        D_{\xi}^{\gamma}c^{\peq{\{\gamma\}}}_x(\xi) 
\end{equation*}
and
\begin{equation*} \label{ccinfprp}
   c^{\peq{\{\gamma\}}}_x(\xi) 
    := 
   c^{\peq{\{\gamma\}}}(x,\xi) \in \ci(\rn \times \rn),
\end{equation*}
and moreover
\begin{equation*} \label{cgamdf}
  c^{\peq{\{\gamma\}}}_x 
    \mbox{ is a smooth function of } x, 
       \mbox{ with values in } \frakX^{m}.   
\end{equation*}
Then by hypothesis \ref{hypo-cuatro}) of this theorem, 
$A^{\peq{\{\gamma\}}} (x) \in 
\Op{\widetilde{\frakX}^{m - \tau \|\gamma\|}}$, so that (\ref{Aastasy}) 
is an asymptotic expansion within the $\Op{\widetilde{\frakX}^{\mu}}$ spaces.
 \end{enumerate}
\end{teo}
\begin{proof}
The proof can be found in \cite{taylor-part-i}.
\end{proof}

\subsection{\sc pseudodifferential calculus on a homogeneous group G} 
\ \\[2mm]
Drawing from 
results of previous sections, 
the following theorem shows that in the case of  a homogeneous group $G$
the family of spaces
  $\left\{\mcall{j}(G)  \right\}_{j\in \rr}$,
   as in Definition \ref{def-multiplicador},
does in fact generate a pseudodifferential calculus, 
fully analogous to the usual pseudodifferential calculus on Euclidean space,
insofar as products and adjoints are concerned. 
  Explicitly, elements of $\Op{\widetilde{\mcal}^j}$
   are operators of the form $\frakA$, where
   \[ 
   \left[ \frakA f \right] (y) = [\check{a}_y \ast f](y), 
   \]
   where $a(x,\xi)$ satisfies, for any compact set $B \subseteq G$,
  the estimates
   \[    
   \Bigabs{D^\beta_x D_\xi^\alpha a(x,\xi)}
   \leq
   C_{\alpha,\beta,B} \big( 1+|\xi|)^{j - |\alpha|}
   \qquad \forall x\in B,\ \xi\in G.  
   \]
 \begin{teo} 
  Suppose $G$ is a homogeneous group, 
  with weights $a_1,$$\ldots$$,a_n$.
   Then  the family of multiplier spaces 
   $\displaystyle \left\{ \mcall{j}(G)\right\}_{j\in\rr}$
   satisfies the following properties
 \begin{enumerate}[\hspace*{1em} a) \  ]
  \item
    $\left\{ \mcall{j}(G) \right\}_{j\in\rbb}  $ is a nested family
 of Fr\'echet spaces.
  \item  If $j\geq 0$ then
        $\displaystyle \mcall{j}(G) \subset \sro{\#}$
         for $\rho=\frac{a_1}{a_n}$.
  \item  If $j< 0$ then
        $\displaystyle \mcall{j}(G) \subset \sroo{m\sigma}{\#}$
         for $\rho=\frac{a_1}{a_n}$,
         and $\sigma=\frac{1}{a_n}$.
  \item
     $ \mcalcopp{j_1}(G) * \mcalcopp{j_2}(G)  \subseteq
      \mcalcopp{j_1+j_2} (G)$
     and the product is continuous.
\item
If $J \in \mcall{j}(G)$, and $\alpha$ is a multiindex, then
  $D^\alpha J\in \mcall{j-\abs{\alpha}}(G)$.  Consequently
  $D^\alpha J\in \mcall{j-\|\alpha\|}(G)$.  
  \item
  Let $J_i \in \mcall{j-i}(G)$, for $i=0,1,2, \ldots$.
 Then there exists a $J\in \mcall{j}(G)$ such that,
 for any $M$, if $N$ is sufficiently large
\[ \left(J-\sum_{i=0}^N J_i\right) \in S^{-M}_{\rho,\#}
     \qquad
         \mbox{for  }  \rho=\frac{a_1}{a_n}.\]
  \item
If $J\in \mcall{j}(G)$ then $\bar{J}\in \mcall{j} (G)$.
 \end{enumerate}
And therefore 
we have 
on $G$ 
a pseudodifferential calculus including products and adjoints,
analogous to 
H\"or\-man\-der's  
$S_{1,0}$-pseudodifferential calculus on $\rn$.  
In fact (by e)), in the situation of
(\ref{ABasymp}), $A^{\peq{[\gamma]}} (x) B_\peq{{[\gamma]}}(x) \in 
\Op{\widetilde{{\mathcal M}}^{m_1 + m_2 - |\gamma|}}$, 
 and in the situation of 
(\ref{Aastasy}), $A^{\peq{\{\gamma\}}} (x) \in 
\Op{\widetilde{{\mathcal M}}^{m - |\gamma|}}$.
\end{teo}
\begin{proof} 
\begin{enumerate}[a)]
      \item
 This is clear. We also note that if  $j_1\leq j_2$, 
 the inclusion map $\mcal^{j_1}(G) \subseteq \mcal^{j_2}(G) $ is continuous.
\addtolength{\leftmargin}{-2em}
\addtolength{\leftmargini}{-2em}
\addtolength{\leftmarginii}{-2em}
\item 
  If $J\in \mcall{j} (G) $, and $\alpha \in (\z^+)^n$ then 
there exists a positive constant $C_\alpha$
such that for all $\xi$
\begin{equation*} 
    \begin{split}
   \bigabs{\partial^\alpha J(\xi) }
       & \leq C_\alpha  (1+\abs{\xi})^{ j-\abs{\alpha} }\\          
       & \leq C_\alpha  
            \frac{(1+\abs{\xi})^j\phantom{xxx} }
               { (1+\abs{\xi})^{ \norm{\alpha} a_1}  }        
        \\
       & \leq C_\alpha'  
            \frac{(1+\abs{\xi})^j\phantom{} }
                {\langle \xi \rangle^{\tfrac{a_1}{a_n} \norm{\alpha}}}\\
       & \leq C_\alpha''  
              {\langle \xi \rangle^{j-\tfrac{a_1}{a_n} \norm{\alpha}}   }\\
      \end{split}
  \end{equation*}
where for the second and third inequalities we have made use 
of Lemma \ref{lema-16}.
Then      $J\in S^j_{\left(\tfrac{a_1}{a_n}\right),\# }$.
 Therefore 
     $\mcall{j}(G) \subseteq S^j_\peq{\rho,\#}$
with $\rho = \frac{a_1}{a_n} \in (0,1]$.
\item 
  If $J\in \mcall{j}(G)$, and $\alpha \in  (\z^+)^n$,
 then there 
  exists a positive constant $C_\alpha$ such that
 for all $\xi$
 \begin{equation*} 
    \begin{split}
        \bigabs{\partial^\alpha J(\xi) }
       & \leq C_\alpha  (1+\abs{\xi})^{ j-\abs{\alpha} }               
       \\ 
      & \leq C_\alpha'  
            \frac{(1+\abs{\xi})^j   }
               { \langle \xi \rangle^{\norm{\alpha}a_1} }
       \\ 
       & \leq C_\alpha''  
           \langle \xi \rangle^\peq{  {\tfrac{1}{a_n}j- 
                                     \frac{a_1}{a_n}\norm{\alpha} 
                                   } }
      \end{split}
  \end{equation*}
where for the second and last inequalities we have made use
of Lemma \ref{lema-16}.
Then
$J\in S_{ \left(\tfrac{a_1}{a_n} \right),\# }^{ j \tfrac{1}{a_n}}$.
Therefore $\mcall{j}(G) \subseteq S^\peq{{j\sigma}}_{\rho,\# }$,
  with
  $\sigma =\frac{1}{a_n}$, and $\rho =\frac{a_1}{a_n}$.
\item 
 In the previous section we established that 
$ \mcalcopp{j_1}(G) * \mcalcopp{j_2}(G)  \subseteq
 \mcalcopp{j_1+j_2} (G)$
for all $j_1, j_2 \in \rbb$. 

To prove the continuity of the product 
we define the following bilinear map.
For a fixed pair $j_1,j_2\in\rbb$,  we set
\begin{equation*}\label{} 
  \xymatrix{
  T\, :\,  \mcall{j_1}(G) \times \mcall{j_2}(G) 
           \ar[r] & \mcall{j_1+j_2} (G)
           }\\
    \end{equation*}\vspace*{-0.4cm}
\begin{equation*}\label{} 
  \xymatrix{
\ \ \ \ \ \ \ \ \ \ \ \ \ \ \ \ 
 ( J_1 ,J_2  ) \ar@{|-{>}}[rr] && {(\copete{J_1}*\copete{J_2})^\sombrero }
           }
    \end{equation*}
and will show that $T$ is continuous.

We consider the following mappings, where ``double'' arrows are used
 to denote mapping which are sequentially continuous
 or separately sequentially continuous
 (for elementary reasons that we will be discussed soon)
\begin{equation*}\label{} 
  \xymatrix{
      &  \scal'(G)   &         \\
  \mcall{j_1}(G) \times \mcall{j_2}(G)\ar@{=>}[ur]^{\widetilde{T} }
     \ar@{->}[r]_-{T}
        &  \mcall{j_1+j_2} (G)
      \ar@{=>}[u]_{i}
           } 
    \end{equation*}
 Here 
$i$ denotes the inclusion.  
$\widetilde{T}$ is 
the same map as $T$ but it is thought 
 as a map from $\mcall{j_1}(G) \times \mcall{j_2}(G)$
 to $\scal'(G)$.
 Notice that the spaces, 
     $ \mcall{j_1+j_2}(G)$,
     and 
     $\mcall{j_1}(G) \times \mcall{j_2}(G)$,
    are Fr\'echet spaces.

To prove the continuity of $T$
 it suffices to prove that 
  the diagram is accurate: 
  $i$~is sequentially continuous,
 and 
$\widetilde{T}$ is separately sequentially continuous. 
 Then by the Closed Graph Theorem, 
 we will know that $T$ is separately sequentially continuous, hence continuous.
 It is clear that $i$ is sequentially continuous. 

 Now we want to prove that $\widetilde{T}$ is separately sequentially
 continuous, which means that we want to show that 
for each fixed 
$J_2 \in   \mcall{j_2}(G)$
and
$J_1 \in   \mcall{j_1}(G)$
the operators
   \begin{equation*}\label{} 
   \xymatrix@1{
      \mcall{j_1}(G) \ 
            \ar@{->}[rr]^{T_{J_2}}
           & &  \scal'(G)
      & &\mcall{j_2}(G) \ 
            \ar@{->}[rr]^{T_{J_1}}
           & &  \scal'(G)
              } 
    \end{equation*}
 \vspace*{-5mm}
 \begin{equation*}\label{} 
   \ \ \ \ \ \ \  \  \xymatrix@1{
         \  \ \ J_1 \ \ \ \
                 \ar@{|->}[rr]
       &   &\big(\copete{J}_1 * \copete{J}_2 \big)^{\sombrero}
      & &   \  \ \ J_2 \ \ \ \
                 \ar@{|->}[rr]
       &   &\big(\copete{J}_1 * \copete{J}_2 \big)^{\sombrero}
              } 
    \end{equation*}
are sequentially continuous.

In order to show that the mapping $T_{J_2} $ is sequentially continuous,
 we shall see that it can be written as the composition of
sequentially continuous mappings $F,H$ defined as follows,
and the Fourier transform. 
   \begin{equation*}\label{} 
   \begin{split}
    \xymatrix@1{
        {\ \mcall{j_1}(G)  \ }
            \ar@{->}[rr]^{F}
           & & {\ \ecal'\oplus \scal \  }
            \ar@{->}[rr]^{H}
            & & {\ \scal'(G) \ }
            \ar@{->}[rr]^{\sombrero}
            & &  {\hspace*{3mm} \scal'(G) \ }
              } 
          \hspace*{9mm} 
              & \\ 
  \xymatrix@1{
           { \scriptstyle 
           J_1  
           }
          \ar@{|->}[rr]
            &&   
           { \scriptstyle 
              \big(\zeta 
                   \copete{J}_1 ,
                   (1-\zeta)  \copete{J}_1  
                  \big)
           }
            \ar@{|->}[r]
            &
           { \scriptstyle 
            \big( \zeta \copete{J}_1\big) \ast \copete{J}_2 
                 + 
                 \left[ (1-\zeta) \copete{J}_1  \right] \ast  \copete{J}_2 
            }
            \ar@{|->}[r]
           & 
             { \scriptstyle 
             \Big[ 
              \big( \zeta \copete{J}_1\big) \ast \copete{J}_2
                 +
                 \left[(1-\zeta) \copete{J}_1 \right] \ast  \copete{J}_2
            \Big]^{\sombrero} 
             }
                 } 
   \end{split}
    \end{equation*}
where $\zeta \in \cisc(G)$, and $\zeta = 1 $ in a neighborhood
of zero.

  $H$ is sequentially continuous, 
  since convolution with $\copete{J}_2 \in \scal'$ is 
  continuous from $\ecal'$ to $\scal'$, or from $\scal'$ to $\scal'$.

Since $\sombrero$ is 
the Fourier transform from $\scal'$ to $\scal'$, 
it is continuous.

Finally in order to show that $F$ is continuous, for a fixed
 $\zeta  \in \cisc(G)$, we shall show that 
 the mappings
   \begin{equation*}\label{} 
    \begin{split}
   \xymatrix@1{
      \mcall{j_1}(G) \ 
            \ar@{->}[rr]^{R_1}
           & &  { \hspace*{2mm}\ecal' }
      & &\mcall{j_1}(G) \ 
            \ar@{->}[rr]^{R_2}
           & &  { \hspace*{2mm} \scal  }
              }\hspace*{7mm} 
\\[-1mm]
\hspace*{4mm} \xymatrix@1{
         \  \ \ J_1 \ \ \ \
                 \ar@{|->}[rr]
       &   &\hspace*{3mm}\zeta  \copete{J}_1 
      & &   \  \ \ J_1 \ \ \ \
                 \ar@{|->}[rr]
       &   &\hspace*{3mm}\big( 1- \zeta \big) \copete{J}_1 
              } 
      \end{split}
    \end{equation*}
are continuous.

Consider the following mappings
\begin{equation*}\label{} 
  \xymatrix{
      &  \scal'(G)  &         \\
  \mcall{j_1}(G) \ar@{=>}[ur]^{\widetilde{R}_2}
     \ar@{->}[r]_-{{R}_2}
        &  \scal(G)
      \ar@{=>}[u]_{i}
           } 
    \end{equation*}
The sequential continuity of $\widetilde{R}_2$ follows from 
  $i: \mcall{j_1}(G)\hookrightarrow \scal'(G)$ 
being sequentially continuous. 
 The continuity of $R_2$ now follows from the the Closed Graph Theorem.

The continuity of $R_1$ also follows from the inclusion
$i: \mcall{j_1}(G)\longrightarrow \scal'(G)$ 
being continuous.

The sequential continuity of $T_{J_1}$ follows in a completely
  analogous fashion. Hence $\widetilde{T}$ is separately
 sequentially continuous, and therefore $T$ is continuous.
 \item
  Assume $\alpha$ is any multiindex. 
  Since $J\in \mcall{j}(G)$, then 
$D^\alpha J \in \mcall{j-|\alpha|} (G) \subseteq 
                    \mcall{j-\|\alpha\|}(G)$.
\item
Let $\psi: \rbb^n\longrightarrow[0,1]$ be a $\ci$ function
  such that $\psi =0 $ for $\abs{\xi} \leq 1$, and $\psi =1$ for 
  $\abs{\xi}\geq 2.$
  Let $\{t_i\}_i$ be a positive decreasing sequence such that
  all $t_i \leq 1$ and 
  $\displaystyle\lim_{i\to 0^+}t_i=0.$
 These $t_i$ will be specified later.
 Define 
      \[ 
     J(\xi) = \sum_{i=0}^\infty \psi \big( \delta_{t_i} (\xi) \big) J_i(\xi) 
    \]
  Since $t_i\longrightarrow 0 $  as $i\longrightarrow 0^+$, then for 
  any fixed $\xi$, $\psi(\delta_{t_i}(\xi) ) = 0 $ for all,
  except for a finite number of $i$; so there are only finitely 
  many non zero terms in the previous sum.  Consequently this sum
 is well defined, and it follows that $J\in \ci (\rbb^n)$.
 In order to prove part f), it suffices to show that 
$J-\sum_{i=0}^N J_i\in \mcall{j-(N+1)} (G) $
  for any $N\in \z^+ $, since  then 
 by parts (b) and (c)  
we shall have 
\[   J- \sum_{i=0}^N J_i \in \mcall{j-(N+1)} (G) \subseteq
    S_{\left(\tfrac{a_1}{a_n}\right),\#}^\peq{ [j-(N+1)] \sigma}  
    \qquad \mbox{ for }  \sigma=\frac{1}{a_n}. \]
     Hence if we choose $N$ sufficiently large so that 
    $[j- (N+1)]\sigma < -M $
\[ J - \sum_{i=0}^N J_i 
    \in 
  S^\peq{-M}_{\left(\tfrac{a_1}{a_n}\right),\#} \]
 as desired.

 For $|\beta|\neq 0,$
       $\big(\partial^\beta \psi \big) \big( \delta_t (\xi)\big) =0 $,
       when $\abs{\delta_t(\xi) } = t \abs{\xi} \leq 1$,
   or when $\abs{\delta_t(\xi)} = {t \abs{\xi}\geq 2}$.
Thus for $\abs{\beta}\neq 0$, \  
 $\big(\partial^\beta \psi \big) \big(\delta_t (\xi)\big) \neq 0$
implies that $1 < \abs{\delta_t(\xi) } = t \abs{\xi} < 2 $, 
i.e. $t^{-1} < \abs{\xi} < 2 t^{-1}$. 
This implies that if $0 < t \leq 1 $ 
there exists some positive constant $C'$ such
 that $t < C' (1+\abs{\xi})^{-1}$.  Therefore if $0< t \leq 1 $, we have
\[ \bigabs{ \partial^\beta [\psi\left( \delta_t(\xi)\right) ]
           }
    =       t^\abs{\beta}
        \bigabs{ (\partial^\beta  \psi) (\delta_t (\xi) )
               }
    \leq 
     C'' t^\abs{\beta}
    \leq
C_\beta (1+\abs{\xi} )^{-\abs{\beta}}\ .
\]
Therefore $\displaystyle \Big\{ \psi\circ \delta_t \Big\}_{0<t \leq 1}$ 
is a bounded subset of $\mcall{0} (G).$

From Leibniz's rule 
and since $J_i \in \mcall{j- i} (G)$,
  it easily follows that
\[  \Bigabs{\partial^\alpha
       \big[
       \psi (\delta_t (\xi))  J_i(\xi)
       \big]
          }
  \leq 
  C_{i,\gamma} ( 1+\abs{\xi} )^{(j- i) -\abs{\alpha}}\ . \]
In particular
$\displaystyle \Big\{ \big( \psi\circ \delta_t \big) J_i 
                         \Big\}_{0<t \leq 1} \subset \mcall{j- i} (G).$

 We set $\displaystyle C_i = \max \big\{C_{i,\alpha} 
                           \ : \ \abs{\alpha} \leq i
                           \big\} $.
 Since $\psi(\xi) = 0$ if $\abs{\xi} \leq 1 $ 
 we have $\psi \big( \delta_t (\xi) \big) \neq 0$, implies 
 $\abs{\delta_t (\xi) } = t \abs{\xi} > 1 $.
 We select $t_i >0$ such that
 $t_i < t_{i-1}$,  and 
  $C_i t_i \leq 2^{-i}. $ 
 Then $\psi \big( \delta_{t_i}( \xi) \big) = 0$ 
 if $ t_i(1+\abs{\xi} ) \leq  1$.
 Therefore for any multiindex $\alpha$ such that $\abs{\alpha} \leq i$
 we have
\begin{equation}\label{ecuacion-cincuenta-y-dos} 
    \begin{split}
     \Bigabs{
       \partial^\alpha \big[    \psi (\delta_{t_i} (\xi) )
            J_i\big] 
            }           & \leq  
            C_i ( 1+ \abs{\xi}  )^{j-  i -\abs{\alpha}} 
        \\ & <  
         C_i t_i ( 1+ \abs{\xi}  )^{j- i +1-\abs{\alpha}}
        \\ & \leq  
       2^{-i}  ( 1+ \abs{\xi}  )^{j-  i +1 -\abs{\alpha}}
     \end{split}
 \end{equation}
For any multiindex $\beta$ we choose $i_o$ such that $i_o \geq\abs{\beta}$,
and we express $J$ as
  \begin{equation}\label{geller-uno} 
   J= \sum_{i=0}^{i_o} \big(\psi\circ \delta_{t_i} \big) J_i  
     + \sum_{i=i_o + 1 }^\infty \big( \psi\circ \delta_{t_i} \big) J_i 
    \end{equation}
Since $
\sum_{i=0}^{i_o}  \big( \psi\circ \delta_{t_i} \big)J_i$
is a finite sum and $  \big( \psi\circ \delta_{t_i} \big)J_i \in 
 \mcall{j- i} (G) \subseteq \mcall{j} (G) $;
we have  
$\sum_{i=0}^{i_o}  \big( \psi\circ \delta_{t_i} \big) J_i 
\in \mcall{j}(G). $
Therefore there exists a positive constant $C$ such that for~all~$\xi$
  \begin{equation}\label{geller-dos} 
    \left\vert \partial^\beta  \sum_{i=0}^{i_o}
                           \psi\big(   \delta_{t_i}(\xi) \big) J_i (\xi)
    \right\vert 
    \leq C ( 1+ \abs{\xi} )^{j-\abs{\beta}}
    \end{equation}
By (\ref{ecuacion-cincuenta-y-dos}), 
  \begin{equation}\label{geller-tres} 
  \left\vert \partial^\beta  \sum_{i=i_o+1}^\infty
                           \psi\big(   \delta_{t_i}(\xi) \big) J_i (\xi)
    \right\vert
   \leq
  \sum_{i=i_o+1}^\infty  2^{-i}
                 \big(1 + \abs{\xi}\big)^{j- i+ 1-\abs{\beta}}
 \leq  \big(1 + \abs{\xi}\big)^{j-\abs{\beta}}\ . 
  \end{equation}
From (\ref{geller-uno}), (\ref{geller-dos}), and (\ref{geller-tres}),
$J\in \mcall{j}(G)$.

Moreover, for $N\in \z^+ $, writing
\[
 J - \sum_{i=0}^N J_i  
 =  \sum_{i=0}^N
   \left[  \big(\psi\circ  \delta_{t_i} \big) -1
   \right] J_i
 + \sum_{i=N+1}^\infty 
   \big(\psi \circ \delta_{t_i} \big) J_i
    \]
and working in the same fashion as before, we obtain
$$\left(\displaystyle \sum_{i=N+1}^\infty 
       (\psi \circ \delta_{t_i} ) J_i \right)\in 
      \mcall{j- (N+1) }(G).
      $$
On the other hand, since 
 $\psi(\xi)=1$      
         for $\abs{\xi}\geq 2 $, 
 $\psi (\delta_{t_i} (\xi) )- 1 = 0 $ 
         for $\abs{\delta_{t_i} (\xi) }= t_i {\abs{\xi} \geq 2}$.
 So, if $0 \leq i \leq N$, $\psi ( \delta_{t_i} (\xi) ))- 1 = 0 $
 for $\abs{\xi} \geq 2 t^{-1}_N.$
 Then
     $
     \sum_{i=0}^N \big[  (\psi \circ \delta_{t_i}) -1 \big] J_i 
     \in \mcall{-\infty}(G) = \scal(G)
     $.
Consequently for any $N$
 \[ 
  \left(
       J - \sum_{i=0}^N J_i 
  \right)  \in \mcall{j- (N+1)} (G)
    \]
i.e. $J \sim \sum_{i=0}^{\infty} J_i$.
\item This is evident.
\qed
\end{enumerate}
\cancelqed
\end{proof}


\end{document}